\input amstex.tex
\input amsppt.sty

\magnification 1000
\pageheight{9truein} \pagewidth{6.5 truein}

\def\bi{\it}

\nologo\TagsOnRight
\NoBlackBoxes

\newif\ifcontains    
\def\contains#1#2{\def\test##1#1##2##3\endtest{
\ifx##2\ZZZ%
\containsfalse\else\containstrue\fi}\test#2#1\ZZZ
\endtest}

\def\cite#1{%
 \def\commasplit##1,##2\endcommasplit{[$
\bold{##1}$,##2]}%
 \contains,{#1}%
 \ifcontains \commasplit#1\endcommasplit%
 \else [{\bf#1}]\fi}

\def\qed{\ifhmode\unskip\nobreak\fi\quad\hfill
\raise.32em\hbox{\boxed{}}}

\topmatter
\title
Variations of
the boundary geometry \\
of 3--dimensional
hyperbolic convex cores
\endtitle
\author
Francis Bonahon
\endauthor
\affil
University of Southern California
\endaffil
\address
Department of Mathematics, 
University of Southern California, 
Los Angeles CA 90089-1113, U.S.A.
\endaddress
\dedicatory
\hfill Dedicated to D.B.A.~Epstein, on his 60th 
birthday.
\enddedicatory
\email fbonahon\@math.usc.edu \endemail

\thanks
This research was partially supported by N.S.F. 
grants DMS-9001895, DMS-9201466 and DMS-9504282.
\endthanks
\endtopmatter\document

\rightheadtext {3--dimensional hyperbolic convex 
cores}

\nobreak\par\nobreak 	Let $M$ be a hyperbolic 
3--manifold, namely a complete 3--dimensional 
Riemannian manifold of constant curvature $-1$, 
such that the fundamental group $\pi _{1}\left ( 
M\right )$ is finitely generated. A fundamental 
subset of $M$ is its {\bi convex core \/} 
$C_{M}$, defined as the smallest non-empty closed 
convex subset of $M$. The boundary $\partial 
C_{M}$ of this convex core is a surface of finite 
topological type, and its geometry was described 
by W.P. Thurston \cite{Thu}: The surface 
$\partial C_{M}$  is almost everywhere totally 
geodesic, and is bent along a family of disjoint 
geodesics called its {\bi pleating locus \/}. The 
path metric induced by the metric of $M$ is 
hyperbolic, and the way $\partial C_{M}$ is bent 
is completely determined by a certain measured 
geodesic lamination.
\nobreak\par\nobreak 	We want to investigate how 
the geometry of $\partial C_{M}$ varies as we 
deform the metric of $M$. For technical reasons, 
in particular because we do not want the topology 
of $\partial C_{M}$ to change, we choose to 
restrict attention to {\bi quasi-isometric 
deformations \/} of $M$, namely hyperbolic 
manifolds $M'$ for which there exists a 
diffeomorphism $M\rightarrow M'$ whose 
differential is uniformly bounded. In the 
language of Kleinian groups, a quasi-isometric 
deformation of $M$ is also equivalent to a 
quasi-conformal deformation of its holonomy; see 
\cite{Thu, {\S }10}. This is not a very strong 
restriction. For instance, in the conjecturally 
generic case where $M$ is geometrically finite 
without cusps, every small deformation of the 
metric is quasi-isometric. When $M$ is 
geometrically finite, quasi-isometric 
deformations of the metric coincide with 
deformations of the holonomy $\pi _{1}\left ( 
M\right )\rightarrow \roman {Isom}\left ( {\Bbb 
H}^{3}\right )$ that respect parabolicity 
\cite{Mar}. Also, every holomorphic family of 
hyperbolic manifolds homeomorphic to $M$ consists 
of quasi-isometric deformations \cite{Su2}.
\nobreak\par\nobreak 	Let ${\Cal Q}{\Cal D}\left 
( M\right )$ be the space of quasi-isometric 
deformations of the metric of $M$, where we 
identify two deformations $M\rightarrow M'$ and 
$M\rightarrow M''$ when the corresponding pull 
back metrics on $M$ are isotopic. This space can 
be parametrized by the space of conformal 
structures on the domain of discontinuity of $M$ 
\cite{Ber}\cite{Su1}, and in particular is a 
differentiable manifold of dimension $3 \left 
\vert \chi \left ( \partial C_{M}\right )\right 
\vert -c$, where $\chi \left ( \enskip \right )$ 
denotes the Euler characteristic and where $c$ is 
the number of cusps of $\partial C_{M}$. Given a 
quasi-isometric deformation $M'$, there is a 
homeomorphism between $\partial C_{M}$ and 
$\partial C_{M'}$, well defined up to isotopy. 
Consequently, if we consider the geometry of 
$\partial C_{M'}$, its hyperbolic metric defines 
an element $\mu \left ( M'\right )$ of the 
Teichm{\"u}ller space ${\Cal T}\left ( \partial 
C_{M}\right )$, and its bending measured geodesic 
lamination defines an element $\beta \left ( 
M'\right )$ of the space ${\Cal M}{\Cal L}\left ( 
\partial C_{M}\right )$ of compact measured 
geodesic laminations on $\partial C_{M}$; see 
\cite{Thu}\cite{CEG}\cite{EpM} for a definition 
of these notions. 
\nobreak\par\nobreak 	Before going any further, 
we must mention that the definitions have to be 
adapted in the special case where the convex core 
$C_{M}$ is a totally geodesic surface, namely 
when $M$ is Fuchsian or twisted Fuchsian. To keep 
the correspondence between $\partial C_{M}$ and 
the domain of discontinuity of $M$, we define in 
this case $\partial C_{M}$ as the unit normal 
bundle of $C_{M}$ in $M$, namely as the `two 
sides' of $C_{M}$ in $M$, whereas the topological 
boundary of $C_{M}$ is equal to $C_{M}$. With 
this convention, we have as above a prefered (up 
to isotopy) identification between $\partial 
C_{M}$ and $\partial C_{M'}$ for every 
quasi-isometric deformation $M\rightarrow M'$, 
and such a deformation again defines a hyperbolic 
metric $\mu \left ( M'\right )\in {\Cal T}\left ( 
\partial C_{M}\right )$ and a bending measured 
lamination $\beta \left ( M'\right )\in {\Cal 
M}{\Cal L}\left ( \partial C_{M}\right )$. 

\proclaim {Theorem 1 }   For every hyperbolic 
$3$--manifold $M$, the map $\mu :{\Cal Q}{\Cal 
D}\left ( M\right )\rightarrow {\Cal T}\left ( 
\partial C_{M}\right )$, defined by considering 
the hyperbolic metrics of convex core boundaries, 
is continuously differentiable. \endproclaim    

\nobreak\par\nobreak 	A simple example in {\S }6 
shows that the map $\mu $ is not necessarily 
twice differentiable. 
\nobreak\par\nobreak 	To prove a similar 
differentiability property for the map $\beta 
:{\Cal Q}{\Cal D}\left ( M\right )\rightarrow 
{\Cal M}{\Cal L}\left ( \partial C_{M}\right )$, 
we encounter a conceptual difficulty. Indeed, the 
space ${\Cal M}{\Cal L}\left ( \partial 
C_{M}\right )$ does not have a natural 
differentiable structure. On the other hand, it 
has a natural structure of piecewise linear 
manifold of dimension $3 \left \vert \chi \left ( 
\partial C_{M}\right )\right \vert -c$; see for 
instance \cite{Thu} \cite{PeH}. In this context, 
we can use a weak notion of differentiability, 
namely the existence of a tangent map (see {\S }1 
for a definition).

\proclaim {Theorem 2 }   The map $\beta :{\Cal 
Q}{\Cal D}\left ( M\right )\rightarrow {\Cal 
M}{\Cal L}\left ( \partial C_{M}\right )$, 
defined by considering the bending measured 
laminations of convex core boundaries, is 
tangentiable in the sense that it admits a 
tangent map everywhere. \endproclaim    

\nobreak\par\nobreak 	The tangent map of $\beta $ 
plays an important r\^ole in the variation of the 
volume of the convex core $C_{M}$, as one varies 
the hyperbolic metric; see \cite{Bo4}. A 
continuity property for the maps $\mu $ and 
$\beta $ was earlier obtained by L.~Keen and 
C.~Series \cite{KeS}.  
\nobreak\par\nobreak 	The proof of Theorems~1 and 
2 is probably of as much interest as the results 
themselves. Indeed, these two statements are 
proved simultaneously, mixing together the 
differentiable and piecewise linear contexts. In 
particular, the `corners' of the piecewise linear 
structure of ${\Cal M}{\Cal L}\left ( \partial 
C_{M}\right )$ account for the fact that the map 
$\mu $ is not $\roman {C}^{2}$.
\nobreak\par\nobreak 	The proof goes as follows. 
First of all, we can restrict attention to the 
case where $M$ is orientable. Indeed, if 
$\widehat M$ is its orientation covering, the 
spaces ${\Cal Q}{\Cal D}\left ( M\right )$, 
${\Cal T}\left ( \partial C_{M}\right )$ and 
${\Cal M}{\Cal L}\left ( \partial C_{M}\right )$ 
are submanifolds (in the appropriate category) of 
${\Cal Q}{\Cal D}\bigl ( \widehat M\bigr ) $, 
${\Cal T}\left ( \partial C_{\widehat M}\right )$ 
and ${\Cal M}{\Cal L}\left ( \partial C_{\widehat 
M}\right )$, respectively, and the maps $\mu $, 
$\beta $ for $M$ are just the restrictions of the 
corresponding maps for $\widehat M$. 
Consequently, we will henceforth assume that $M$ 
is orientable.
\nobreak\par\nobreak 	Let $S_{1}$, \dots , 
$S_{n}$ be the components of $\partial C_{M}$. 
For each $i$, let ${\Cal R}\left ( S_{i}\right )$ 
denote the space of representations $\pi 
_{1}\left ( S_{i}\right )\rightarrow \roman 
{Isom}^{+}\left ( {\Bbb H}^{3}\right )$ sending 
the fundamental group of each end of $S_{i}$ to a 
parabolic subgroup of $\roman {Isom}^{+}\left ( 
{\Bbb H}^{3}\right )$, where $\roman 
{Isom}^{+}\left ( {\Bbb H}^{3}\right )$ denotes 
the group of orientation-preserving isometries of 
the hyperbolic 3--space ${\Bbb H}^{3}$ and where 
these representations are considered modulo 
conjugation by elements of $\roman 
{Isom}^{+}\left ( {\Bbb H}^{3}\right )$. Let 
${\Cal R}\left ( \partial C_{M}\right )$ denote 
the product $\prod _{i=1}^{n}{\Cal R}\left ( 
S_{i}\right )$. Restricting the holonomy of a 
quasi-isometric deformation to the components of 
$\partial C_{M}$, we get a map $R:{\Cal Q}{\Cal 
D}\left ( M\right )\rightarrow {\Cal R}\left ( 
\partial C_{M}\right )$. The image of $R$ is in 
the non-singular part of ${\Cal R}\left ( 
\partial C_{M}\right )$, and $R$ is 
differentiable; see for instance \cite{CuS}.
\nobreak\par\nobreak 	If we are given a finite 
area hyperbolic metric and a compactly supported 
measured geodesic lamination on the surface 
$S_{i}$, we can always realize these in a unique 
way as the pull back metric and the bending 
measured lamination of a pleated surface $f=\bigl 
( \widetilde  f, \rho \bigr ) $, where $\rho \in 
{\Cal R}\left ( S_{i}\right )$ is not necessarily 
discrete and where $\widetilde  f:\widetilde  
S_{i}\rightarrow {\Bbb H}^{3}$ is a $\rho 
$--equivariant pleated surface from the universal 
covering of $S$ into ${\Bbb H}^{3}$; see 
\cite{EpM}\cite{KaT}\cite{Bo3}. By considering 
the corresponding representations, this defines a 
map $\varphi :{\Cal T}\left ( \partial 
C_{M}\right )\times {\Cal M}{\Cal L}\left ( 
\partial C_{M}\right )\rightarrow {\Cal R}\left ( 
\partial C_{M}\right )$.  Thurston showed that 
$\varphi $ is a local homeomorphism, by 
establishing a correspondence between ${\Cal 
T}\left ( \partial C_{M}\right )\times {\Cal 
M}{\Cal L}\left ( \partial C_{M}\right )$ and the 
space of complex projective structures on 
$\partial C_{M}$; see \cite{KaT}, and see 
\cite{Kap} for a description of the image of 
$\varphi $. In particular, there is a local 
inverse $\varphi ^{-1}$ defined near the point of 
${\Cal R}\left ( \partial C_{M}\right )$ 
corresponding to the original metric of $M$. 
Then, near that metric, the product $\mu \times 
\beta :{\Cal Q}{\Cal D}\left ( M\right 
)\rightarrow {\Cal T}\left ( \partial C_{M}\right 
)\times {\Cal M}{\Cal L}\left ( \partial 
C_{M}\right )$ coincides with the composition 
$\varphi ^{-1}\circ R$.
\nobreak\par\nobreak 	The main technical step in 
the proof of Theorems~1 and 2 is to show that the 
map $\varphi $ is tangentiable, and that its 
tangent map is everywhere injective. This is done 
in {\S }{\S }2--3, by locally comparing $\varphi 
$ to the parametrization of ${\Cal R}\left ( 
\partial C_{M}\right )$ by shear-bend coordinates 
developed in \cite{Bo3}. The crucial technical 
step here is the growth estimate provided by 
Lemma~7. Then, an easy inverse function theorem 
(Lemma~4 in {\S }1) shows that the local inverse 
$\varphi ^{-1}$ is tangentiable. Since $\mu 
\times \beta =\varphi ^{-1}\circ R$ and $R$ is 
differentiable, it follows that $\mu $ and $\beta 
$ are tangentiable. In addition, the proof gives 
that the tangent map of $\mu $ is linear, so that 
$\mu $ is differentiable in the usual sense. 
Continuity properties for the differential of 
$\mu $ follow from the computation of this 
differential, and are proved in {\S }5.
\nobreak\par\nobreak 	As a by-product of the 
proof, we obtain the following result for the 
space of complex projective structures on a 
connected surface $S$ of finite type (withouth 
boundary). A {\bi complex projective structure 
\/} on $S$ is an atlas modelling $S$ over open 
subsets of the complex projective line ${\Bbb 
C}{\Bbb P}^{1}$, where all changes of charts 
extend to elements of the projective group 
$\roman {PSL}_{2}\left ( {\Bbb C}\right )$ and 
where the atlas is maximal for this property. Let 
${\Cal P}\left ( S\right )$ be the space of 
isotopy classes of complex projective structures 
on $S$ which are of cusp type near the ends of 
$S$. When $\chi \left ( S\right )<0$, Thurston 
defined a homeomorphism $\psi :{\Cal T}\left ( 
S\right )\times {\Cal M}{\Cal L}\left ( S\right 
)\rightarrow {\Cal P}\left ( S\right )$, by 
associating to each complex projective structure 
a locally convex pleated surface; see \cite{KaT} 
for an exposition. Because geometric structures 
are locally parametrized by deformations of their 
monodromy, the monodromy map ${\Cal P}\left ( 
S\right )\rightarrow {\Cal R}\left ( S\right )$ 
is a local diffeomorphism. Our proof that 
$\varphi :{\Cal T}\left ( \partial C_{M}\right 
)\times {\Cal M}{\Cal L}\left ( \partial 
C_{M}\right )\rightarrow {\Cal R}\left ( \partial 
C_{M}\right )$ and its local inverses are 
tangentiable immediately gives:
\proclaim {Theorem 3 }   The Thurston 
homeomorphism $\psi :{\Cal T}\left ( S\right 
)\times {\Cal M}{\Cal L}\left ( S\right 
)\rightarrow {\Cal P}\left ( S\right )$ and its 
inverse are tangentiable. \endproclaim    

\nobreak\par\nobreak 	Again, if we compose $\psi 
^{-1}$ with the projection ${\Cal T}\left ( 
S\right )\times {\Cal M}{\Cal L}\left ( S\right 
)\rightarrow {\Cal T}\left ( S\right )$, the map 
${\Cal P}\left ( S\right )\rightarrow {\Cal 
T}\left ( S\right )$ so obtained is $\roman 
{C}^{1}$ but not $\roman {C}^{2}$. 

\remark {Acknowledgements}   Parts of this 
article were written while the author was 
visiting the University of California at 
Berkeley, the Centre \'Emile Borel and the 
Institut des Hautes \'Etudes Scientifiques. He 
would like to thank these institutions for their 
beneficial hospitality. He is also grateful to 
Mike Wolf for pointing out to him the reference 
\cite{KaT}. \endremark

\head {\S }1. Tangent maps\endhead    

\nobreak\par\nobreak 	Given a map $\varphi 
:U\rightarrow {\Bbb R}^{p}$ defined on an open 
subset $U$ of ${\Bbb R}^{n}$, its {\bi tangent 
map \/} at $x\in U$ is, if it exists, the map 
$T_{x}\varphi :{\Bbb R}^{n}\rightarrow {\Bbb 
R}^{p}$ such that one of the following equivalent 
conditions hold:
	\itemitem {(i)}   $T_{x}\varphi \left ( v\right 
)=\lim_{t\rightarrow 0^{+}} \left ( \varphi \left 
( x+tv\right )-\varphi \left ( x\right )\right 
)/t$, uniformly in $v$ on compact subsets of 
${\Bbb R}^{n}$;
	\itemitem {(ii)}   for every curve $\gamma 
:\left [0,\varepsilon \right [\rightarrow U$ with 
$\gamma \left ( 0\right )=x$ and $\gamma '\left ( 
0\right )=v$, then $T_{x}\varphi \left ( v\right 
)=\left ( \varphi \circ \gamma \right )'\left ( 
0\right )$.
	\itemitem {(iii)}   for every sequence of points 
$x_{n}\in {\Bbb R}^{n}$ and numbers $t_{n}>0$ 
such that $\lim_{n\rightarrow \infty } t_{n}=0$ 
and $\lim_{n\rightarrow \infty } \left ( 
x_{n}-x\right )/t_{n}=v$, then $T_{x}\varphi 
\left ( v\right )=\lim_{n\rightarrow \infty } 
\left ( \varphi \left ( x_{n}\right )-\varphi 
\left ( x\right )\right )/t_{n}$ (a discrete 
version of \rom{(ii)}). 
\nobreak\par\nobreak 	The equivalence of these 
three conditions is an easy exercise. The tangent 
map $T_{x}\varphi $ is positive homogeneous of 
degree 1 (namely $T_{x}\varphi \left ( av\right 
)=aT_{x}\varphi \left ( v\right )$ for every 
$v\in {\Bbb R}^{n}$ and $a\geqslant 0$), but not 
necessarily linear. We will say that $\varphi $ 
is {\bi tangentiable \/} if it admits a tangent 
map at each $x\in U$.
\nobreak\par\nobreak 	A {\bi tangentiable 
structure \/} on a topological manifold is a 
maximal atlas where all changes of charts are 
tangentiable. Examples of such tangentiable 
manifolds include differentiable manifolds, 
piecewise linear manifolds and products of these, 
as we will encounter in this paper.  By the usual 
tricks, we can define a space $T_{x}M$ of tangent 
vectors at each point $x$ of a tangentiable 
manifold $M$. This tangent space $T_{x}M$ is not 
necessarily a vector space, although it admits a 
law of multiplication by non-negative numbers. 
There is also a notion of tangentiable map 
between tangentiable manifolds, defined using 
charts, and such a tangentiable map $\varphi 
:M\rightarrow N$ induces a tangent map 
$T_{x}\varphi :T_{x}M\rightarrow T_{\varphi \left 
( x\right )}N$ for every $x\in M$. 
 
\proclaim {Lemma 4 }   Let $\varphi :M\rightarrow 
N$ be a homeomorphism between two tangentiable 
manifolds. Assume that $\varphi $ admits a 
tangent map at $x\in M$, and that this tangent 
map $T_{x}\varphi :T_{x}M\rightarrow T_{\varphi 
\left ( x\right )}N$ is injective. Then, the 
inverse $\varphi ^{-1}$ admits a tangent map at 
$\varphi \left ( x\right )$, and $T_{\varphi 
\left ( x\right )}\varphi ^{-1}=\left ( 
T_{x}\varphi \right )^{-1}$.       \endproclaim   
\demo {Proof}   Because $\varphi $ is a 
homeomorphism, $T_{x}\varphi $ is surjective by a 
degree argument. The fact that $T_{\varphi \left 
( x\right )}\varphi ^{-1}=\left ( T_{x}\varphi 
\right )^{-1}$ easily follows by taking 
appropriate subsequences in 
Definition~\rom{(iii)} of tangentiability.  
\qed\enddemo

\head {\S }2. Proof that $\varphi :{\Cal T}\left 
( S\right )\times {\Cal M}{\Cal L}\left ( S\right 
)\rightarrow {\Cal R}\left ( S\right )$ is 
tangentiable\endhead    

\nobreak\par\nobreak 	Let $S$ be a connected 
oriented surface of finite type and negative 
Euler characteristic. Given a finite volume 
hyperbolic metric $m$ and a compactly supported 
measured geodesic lamination $b$ on $S$, there is 
a unique locally convex pleated surface $f=\bigl 
( \widetilde  f, \rho \bigr ) $ whose pull back 
metric is equal to $m$ and whose bending measured 
lamination is equal to $b$; see 
\cite{EpM}\cite{KaT}\cite{Bo3}. This defines a 
map $\varphi :{\Cal T}\left ( S\right )\times 
{\Cal M}{\Cal L}\left ( S\right )\rightarrow 
{\Cal R}\left ( S\right )$. This bending map 
$\varphi $ is also the composition of the 
Thurston parametrization $\psi :{\Cal T}\left ( 
S\right )\times {\Cal M}{\Cal L}\left ( S\right 
)\rightarrow {\Cal P}\left ( S\right )$ with the 
holonomy map ${\Cal P}\left ( S\right 
)\rightarrow {\Cal R}\left ( S\right )$. Because 
$\psi $ and the monodromy map ${\Cal P}\left ( 
S\right )\rightarrow {\Cal R}\left ( S\right )$ 
are local homeomorphisms, so is $\varphi $.
\nobreak\par\nobreak 	In \cite{Bo3}, we 
developped another local parametrization of 
${\Cal R}\left ( S\right )$ which similarly uses 
pleated surfaces. Fix a compact geodesic 
lamination $\lambda $ on $S$. If $f=\bigl ( 
\widetilde  f, \rho \bigr ) $ is a pleated 
surface with pleating locus $\lambda $, the 
amount by which $f$ bends along $\lambda $ is 
measured by a transverse finitely additive 
measure for $\lambda $, valued in ${\Bbb R}/2\pi 
{\Bbb Z}$. We call such a transverse finitely 
additive measure an ${\Bbb R}/2\pi {\Bbb 
Z}$--valued {\bi transverse cocycle \/} for 
$\lambda $. In general, this {\bi bending 
transverse cocycle \/} is not a (countably 
additive) transverse measure, unless the pleated 
surface is {\bi locally convex \/}, namely always 
bends in the same direction. Let ${\Cal H}\left ( 
\lambda ;{\Bbb R}/2\pi {\Bbb Z}\right )$ denote 
the space of all ${\Bbb R}/2\pi {\Bbb Z}$--valued 
transverse cocycles for $\lambda $. 
\nobreak\par\nobreak 	Given $m\in {\Cal T}\left ( 
S\right )$ and $b\in {\Cal H}\left ( \lambda 
;{\Bbb R}/2\pi {\Bbb Z}\right )$, there is a 
unique pleated surface $f=\bigl ( \widetilde  f, 
\rho \bigr ) $ with pleating locus $\lambda $, 
pull back metric $m$ and bending transverse 
cocycle $b$. This defines a differentiable map 
$\varphi _{\lambda }:{\Cal T}\left ( S\right 
)\times {\Cal H}\left ( \lambda ;{\Bbb R}/2\pi 
{\Bbb Z}\right )\rightarrow {\Cal R}\left ( 
S\right )$. If, in addition, $\lambda $ is {\bi 
maximal \/} among all compact geodesic 
laminations (this is equivalent to say that each 
component of $S-\lambda $ is, either an infinite 
triangle, or an annulus leading to a cusp and 
with exactly one spike in its boundary), then 
$\varphi _{\lambda }$ is a local diffeomorphism; 
see \cite{Bo3}.
\nobreak\par\nobreak 	Transverse cocycles 
occurred in a different context in \cite{Bo2}. 
The piecewise linear structure of ${\Cal M}{\Cal 
L}\left ( S\right )$ defines a space of tangent 
vectors at each of its points, as in {\S }1. In 
\cite{Bo2}, we gave an interpretation of these 
combinatorial tangent vectors at $a\in {\Cal 
M}{\Cal L}\left ( S\right )$ as geodesic 
laminations containing the support of $a$ and 
endowed with transverse ${\Bbb R}$--valued 
cocycles. In this context, Proposition~5 below 
connects the infinitesimal properties of the maps 
$\varphi :{\Cal T}\left ( S\right )\times {\Cal 
M}{\Cal L}\left ( S\right )\rightarrow {\Cal 
R}\left ( S\right )$ and $\varphi _{\lambda 
}:{\Cal T}\left ( S\right )\times {\Cal H}\left ( 
\lambda ;{\Bbb R}/2\pi {\Bbb Z}\right 
)\rightarrow {\Cal R}\left ( S\right )$.
\nobreak\par\nobreak 	Before stating this result, 
it is convenient to introduce the following 
notation. We will often have to consider the 
right derivatives at $t=0$ of various quantities 
$a_{t}$ defined for $t\in \left [0,\varepsilon 
\right [$,  with $\varepsilon >0$. We will denote 
such a derivative $da_{t}/dt^{+}_{ \vert t=0}$ by 
$\dot a_{0}$. 

\proclaim {Proposition 5 }   Let the 
$1$--parameter families $m_{t}\in {\Cal T}\left ( 
S\right )$ and $b_{t}\in {\Cal M}{\Cal L}\left ( 
S\right )$, $t\in \left [0,\varepsilon \right [$, 
admit tangent vectors $\dot m_{0}$ and $\dot 
b_{0}$ at $t=0$, respectively, and let $\rho 
_{t}=\varphi \left ( m_{t},b_{t}\right )\in {\Cal 
R}\left ( S\right )$. Interpret $\dot b_{0}$ as a 
geodesic lamination with a transverse  ${\Bbb 
R}$--valued cocycle, and choose a maximal 
geodesic lamination $\lambda $ which contains the 
supports of $b_{0}$ and $\dot b_{0}$. In 
particular, $b_{0}$ and $\dot b_{0}$ can both be 
considered as elements of ${\Cal H}\left ( 
\lambda ;{\Bbb R}\right )$, and $\rho 
_{0}=\varphi _{\lambda }\left ( m_{0},\bar 
b_{0}\right )$ where $\bar b_{0}\in {\Cal H}\left 
( \lambda ;{\Bbb R}/2\pi {\Bbb Z}\right )$ is the 
reduction of $b_{0}$ modulo $2\pi $. Then, the 
family $\rho _{t}$ admits a tangent vector $\dot 
\rho _{0}$ at $t=0^{+}$ and $\dot \rho 
_{0}=T_{\left ( m_{0}, \bar b_{0}\right )} 
\varphi _{\lambda }\bigl ( \dot m_{0}, \dot 
b_{0}\bigr ) $.  \endproclaim    

\nobreak\par\nobreak 	The tangent space 
$T_{b_{0}}{\Cal M}{\Cal L}\left ( S\right )$ 
admits a decomposition into linear faces. Each 
face is associated to a geodesic lamination 
$\lambda $ containing the support of $b_{0}$, and 
the tangent vectors in this face correspond to 
(some) transverse cocycles in ${\Cal H}\left ( 
\lambda ;{\Bbb R}\right )$; see \cite{Bo2, {\S 
}5}. Proposition~5 immediately implies the 
following corollary.

\proclaim {Corollary~6 }   The map $\varphi 
:{\Cal T}\left ( S\right )\times {\Cal M}{\Cal 
L}\left ( S\right )\rightarrow {\Cal R}\left ( 
S\right )$ is tangentiable at each $\left ( 
m_{0},b_{0}\right )$. In addition, if $\lambda $ 
is a maximal geodesic lamination containing the 
support of $b_{0}$ and if $\bar b_{0}\in {\Cal 
H}\left ( \lambda ;{\Bbb R}/2\pi {\Bbb Z}\right 
)$ denotes the reduction of $b_{0}$ modulo $2\pi 
$, the tangent map $T_{\left ( m_{0},b_{0}\right 
)}\varphi $ coincides with $T_{\left ( m_{0}, 
\bar b_{0}\right )} \varphi _{\lambda }$ on the 
product of $T_{m_{0}}{\Cal T}\left ( S\right )$ 
and of the face of $T_{b_{0}}{\Cal M}{\Cal 
L}\left ( S\right )$ associated to $\lambda $.  
\qed  \endproclaim    

\demo {Proof of Proposition ~5}   Consider the 
transverse cocycle $b_{t}'=b_{0}+t\dot b_{0}\in 
{\Cal H}\left ( \lambda ;{\Bbb R}\right )$ and 
its reduction $\bar b_{t}'\in {\Cal H}\left ( 
\lambda ;{\Bbb R}/2\pi {\Bbb Z}\right )$ modulo 
$2\pi $. Let $\rho _{t}'=\varphi _{\lambda }\bigl 
( m_{t},\bar b_{t}'\bigr ) \in {\Cal R}\left ( 
S\right )$. Because $\dot b_{0}'=\dot b_{0}$ and 
because $\varphi _{\lambda }$ is a differentiable 
map, the curve $t\mapsto \rho _{t}'$ admits a 
tangent vector $\dot \rho _{0}'=T_{\left ( m_{0}, 
\bar b_{0}\right )} \varphi _{\lambda }\bigl ( 
\dot m_{0},\dot b_{0}\bigr ) $ at $t=0$. We will 
compare the two curves $t\mapsto \rho _{t}$ and 
$t\mapsto \rho _{t}'$ in ${\Cal R}\left ( S\right 
)$, and show that they are tangent at $t=0$.
\nobreak\par\nobreak 	We first make the 
additional assumption that, for the Hausdorff 
topology, the geodesic lamination $\lambda _{t}$ 
underlying $b_{t}$ converges to some 
sublamination of $\lambda $ as $t$ tends to 
$0^{+}$. We will later indicate how to obtain the 
general case from this one.
\nobreak\par\nobreak 	Let $f_{t}=\bigl ( 
\widetilde  f_{t}, \rho _{t}\bigr ) $ be the 
locally convex pleated surface with pull back 
metric $m_{t}$ and bending measured lamination 
$b_{t}$. Similarly, let $f_{t}'=\bigl ( 
\widetilde  f_{t}', \rho _{t}'\bigr ) $ be the 
pleated surface pleated along $\lambda $ with 
pull back metric $m_{t}$ and bending transverse 
cocycle $b_{t}'$. In the universal covering 
$\widetilde  S$, consider the preimage 
$\widetilde  \lambda $ of $\lambda $.
\nobreak\par\nobreak 	So far, the metric $m_{t}$ 
was defined only up to isotopy of $S$, and 
$\widetilde  f_{t}$, $\rho _{t}$, $\widetilde  
f_{t}'$ and $\rho _{t}'$ were only defined up to 
conjugacy by isometries of ${\Bbb H}^{3}$. We can 
normalize these so that the metric $m_{t}$ 
$\roman {C}^{\infty }$--converges to $m_{0}$ and 
so that, for a choice of a base point $\widetilde 
 x_{0}\in \widetilde  S-\widetilde  \lambda $ and 
of a base frame at $\widetilde  x_{0}$, 
$\widetilde  f_{t}$ and $\widetilde  f_{t}'$ 
coincide with $\widetilde  f_{0}$ at these base 
point and frame. 
\nobreak\par\nobreak 	To show that the two curves 
$t\mapsto \rho _{t}$ and $t\mapsto \rho _{t}'$ 
are tangent at $t=0$ in ${\Cal R}\left ( S\right 
)$, it then suffices to show that, for each $\xi 
\in \pi _{1}\left ( S\right )$, the curves 
$t\mapsto \rho _{t}\left ( \xi \right )$ and 
$t\mapsto \rho _{t}'\left ( \xi \right )$ are 
tangent at $t=0$ in $\roman {Isom}^{+}\left ( 
{\Bbb H}^{3}\right )$. For this, we first have to 
remind the reader of the construction of $\bigl ( 
\widetilde  f_{t},\rho _{t}\bigr ) $ and $\bigl ( 
\widetilde  f_{t}',\rho _{t}'\bigr ) $. 
\nobreak\par\nobreak 	We begin with the totally 
geodesic (un-)pleated surface $\bigl ( \widetilde 
 f''_{t},\rho ''_{t}\bigr ) $ with pull back 
metric $m_{t}$ and bending measured lamination 0, 
normalized so that $\widetilde  f''_{t}$ 
coincides with $\widetilde  f_{0}$ at the base 
frame in $\widetilde  S$. To fix ideas, we can 
arrange that $\widetilde  f''_{t}\bigl ( 
\widetilde  S\bigr ) ={\Bbb H}^{2}\subset {\Bbb 
H}^{3}$. Choose as base point $x_{0}$ for the 
fundamental group $\pi _{1}\left ( S\right )=\pi 
_{1}\left ( S;x_{0}\right )$ the image of the 
base point $\widetilde  x_{0}\in \widetilde  S$. 
Let $\widetilde  c$ be the $m_{0}$--geodesic arc 
in $\widetilde  S$ going from $\widetilde  x_{0}$ 
to $\xi \widetilde  x_{0}$, so that the 
projection of $\widetilde  c$ to $S$ represents 
$\xi \in \pi _{1}\left ( S;x_{0}\right )$. Then, 
$\rho _{t}\left ( \xi \right )$ and $\rho 
_{t}'\left ( \xi \right )$ are defined by 
composition of $\rho ''_{t}\left ( \xi \right )$ 
with rotations around certain geodesics of ${\Bbb 
H}^{2}\subset {\Bbb H}^{3}$ that are determined 
by $\xi $, $\lambda $, and $b_{t}$. 
\nobreak\par\nobreak 	Let $U\subset S$ be a train 
track neighborhood carrying $\lambda $, or more 
precisely carrying the $m_{0}$--geodesic 
lamination corresponding to $\lambda $. We can 
choose $U$ sufficiently small so that, if 
$\widetilde  U$ is its preimage in $\widetilde  
S$, each component of $\widetilde  c\cap 
\widetilde  U$ is an arc contained in a single 
edge of $\widetilde  U$. Because of our 
assumption that the $m_{0}$--geodesic lamination 
underlying $b_{t}$ converges to some 
sublamination of $\lambda $, $U$ will also carry 
this lamination for $t$ sufficiently small. 
Finally, since the metric $m_{t}$ converges to 
$m_{0}$, the $m_{t}$--geodesic representative of 
the geodesic lamination underlying $b_{t}$ will 
also be carried by $U$ for $t$ sufficiently 
small. 
\nobreak\par\nobreak 	For $r\geqslant 0$, let 
$\Gamma _{r}$ be the set of all edge paths of 
length $2r+1$ in $\widetilde  U$ that are 
centered on an edge meeting $\widetilde  c$. We 
can partially order the elements of $\Gamma _{r}$ 
from $\widetilde  x_{0}$ to $\xi \widetilde  
x_{0}$ as follows. For two edge paths $\gamma $, 
$\gamma '$ centered at different edges of 
$\widetilde  U$, $\gamma \prec \gamma '$ 
precisely when the central edge of $\gamma $ cuts 
$\widetilde  c$ closer to $\widetilde  x_{0}$ 
than the central edge of $\gamma '$. Two edge 
paths $\gamma $, $\gamma '$ with the same central 
edge $e$ follow a common edge path and diverge at 
1 or 2 switches; then $\gamma \prec \gamma '$ 
precisely when $\gamma $ diverges always on the 
side of $\gamma '$ which contains the point of 
$e\cap \widetilde  c$ that is closest to 
$\widetilde  x_{0}$.  Neither $\gamma \prec 
\gamma '$ nor $\gamma '\prec \gamma $ hold when 
$\gamma $ and $\gamma '$ have the same central 
edge and diverge on opposite sides.  
\nobreak\par\nobreak 	List all the elements of 
$\Gamma _{r}$ as $\gamma _{1}$, $\gamma _{2}$, 
\dots , $\gamma _{p}$ in a way which is 
compatible with the partial order $\prec $, 
namely so that $i<j$ whenever $\gamma _{i}\prec 
\gamma _{j}$. For each $\gamma _{i}$, let 
$g_{i}^{t}$ be the geodesic of ${\Bbb 
H}^{2}\subset {\Bbb H}^{3}$ image under 
$\widetilde  f_{t}'':\widetilde  S\rightarrow 
{\Bbb H}^{2}\subset {\Bbb H}^{3}$ of an 
$m_{t}$--geodesic of $\widetilde  S$ that is 
carried by $\widetilde  U$ and realizes $\gamma 
_{i}$. Such a geodesic may not exist for every 
$\gamma _{i}$, but it will definitely exist if at 
least one of $b_{t}\left ( \gamma _{i}\right )$ 
or $b_{t}'\left ( \gamma _{i}\right )$ is 
non-zero (for instance, a leaf of the 
$m_{t}$--geodesic lamination underlying $b_{t}$ 
if $b_{t}\left ( \gamma _{i}\right )\mathbin{\not 
=}0$, or a leaf of the $m_{t}$--geodesic 
lamination corresponding to $\lambda $ if 
$b_{t}'\left ( \gamma _{i}\right )\mathbin{\not 
=}0$), which is exactly the case in which we need 
it. 
\nobreak\par\nobreak  To each edge path $\gamma $ 
of $\widetilde  U$, the transverse measure of 
$b_{t}$ associates a number $b_{t}\left ( \gamma 
\right )\geqslant 0$, namely the $b_{t}$--mass of 
the set of those geodesics realizing $\gamma $ 
(whether we consider $m_{t}$-- or 
$m_{0}$--geodesics does not matter here because 
the $m_{0}$--geodesic lamination and 
$m_{t}$--geodesic lamination underlying $b_{t}$ 
are both carried by $U$). This $b_{t}\left ( 
\gamma \right )$ is a piecewise linear function 
of $b_{t}\in {\Cal M}{\Cal L}\left ( S\right )$, 
and the fact that $t\mapsto b_{t}$ admits a 
tangent vector at $t=0^{+}$ is equivalent to the 
property that $t\mapsto b_{t}\left ( \gamma 
\right )$ admits a right derivative $\dot 
b_{0}\left ( \gamma \right )$ for every edge path 
$\gamma $. The transverse cocycle $b_{t}'$ 
similarly associates a number $b_{t}'\left ( 
\gamma \right )$ to $\gamma $ which, in our case, 
is equal to $b_{0}\left ( \gamma \right )+t\dot 
b_{0}\left ( \gamma \right )$. See 
\cite{Bo1}\cite{Bo2}. Then,
$$\rho _{t}\left ( \xi \right 
)=\lim_{r\rightarrow \infty } 
R_{g_{1}^{t}}^{b_{t}\left ( \gamma _{1}\right )} 
R_{g_{2}^{t}}^{b_{t}\left ( \gamma _{2}\right )} 
\dots R_{g_{p}^{t}}^{b_{t}\left ( \gamma 
_{p}\right )} \rho ''_{t}\left ( \xi \right )  
\tag 1$$
and
$$\rho _{t}'\left ( \xi \right 
)=\lim_{r\rightarrow \infty } 
R_{g_{1}^{t}}^{b_{t}'\left ( \gamma _{1}\right )} 
R_{g_{2}^{t}}^{b_{t}'\left ( \gamma _{2}\right )} 
\dots R_{g_{p}^{t}}^{b_{t}'\left ( \gamma 
_{p}\right )} \rho ''_{t}\left ( \xi \right ).  
\tag 2$$
where $R_{g}^{b}\in \roman {Isom}^{+}\left ( 
{\Bbb H}^{3}\right )$ denotes the hyperbolic 
rotation of angle $b\in {\Bbb R}/2\pi {\Bbb Z}$ 
around the oriented geodesic $g$, and where the 
$g_{i}^{t}$ are oriented to the left as seen from 
the base point $\widetilde  f_{0}\bigl ( 
\widetilde  x_{0}\bigr ) $ in ${\Bbb H}^{2}$. 
Compare \cite{EpM, {\S }3} for the case of 
transverse measures, and see \cite{Bo3, {\S }5} 
for the more general case of transverse cocycles, 
where the convergence is much more subtle.
\nobreak\par\nobreak 	Identify the isometry group 
$\roman {Isom}^{+}\left ( {\Bbb H}^{3}\right )$ 
to some matrix group, for instance $\roman 
{SO}\left ( 3,1\right )$, and endow the 
corresponding space of matrices with any of the 
classical norms $\left \Vert \enskip \right \Vert 
$ such that $\left \Vert AB\right \Vert \leqslant 
\left \Vert A\right \Vert \left \Vert B\right 
\Vert $.
\nobreak\par\nobreak 	We can write the difference 
$\rho _{t}\left ( \xi \right )-\rho _{t}'\left ( 
\xi \right )$ as
$$\rho _{t}\left ( \xi \right )-\rho _{t}'\left ( 
\xi \right )= \lim_{r\rightarrow \infty } 
A_{r}^{t}-B_{r}^{t} =\lim_{r\rightarrow \infty } 
C_{r}^{t}$$
where 
$$\align
      A_{r}^{t} 	&= R_{g_{1}^{t}}^{b_{t}\left ( 
\gamma _{1}\right )} R_{g_{2}^{t}}^{b_{t}\left ( 
\gamma _{2}\right )} \dots 
R_{g_{p}^{t}}^{b_{t}\left ( \gamma _{p}\right )} 
\rho ''_{t}\left ( \xi \right ) \\
       B_{r}^{t} 	&= R_{g_{1}^{t}}^{b_{t}'\left ( 
\gamma _{1}\right )} R_{g_{2}^{t}}^{b_{t}'\left ( 
\gamma _{2}\right )} \dots 
R_{g_{p}^{t}}^{b_{t}'\left ( \gamma _{p}\right )} 
\rho ''_{t}\left ( \xi \right )
\endalign$$
and $C_{r}^{t}=A_{r}^{t}-B_{r}^{t}$. 
\nobreak\par\nobreak 	The following growth 
estimate is the technical key to the proof of 
Proposition~5.

\proclaim {Lemma 7 }   There is a number $A>0$ 
such that
$$C_{r+1}^{t}-C_{r}^{t}= tO\bigl ( e^{-Ar}\bigl 
\Vert \dot b_{0}\bigr \Vert _{U}\bigr ) $$
and
$$\rho _{t}\left ( \xi \right )-\rho _{t}'\left ( 
\xi \right )=C_{r}^{t} + t O\bigl ( e^{-Ar}\bigl 
\Vert \dot b_{0}\bigr \Vert _{U}\bigr ) $$
where $\bigl \Vert \dot b_{0}\bigr \Vert _{U}$ 
denote the maximum of $\bigl \vert \dot 
b_{0}\left ( e\right )\bigr \vert $ as $e$ ranges 
over all edges of $U$, and where $A$ and the 
constants hidden in the symbols $O\left ( \enskip 
\right )$ are independent of $r$ and $t$. 
\endproclaim    
\demo {Proof of Lemma~7}   List the edge paths of 
$\Gamma _{r+1}$ as $\delta _{1}$, \dots , $\delta 
_{q}$, where the indexing is chosen to be 
compatible with the partial order $\prec $. There 
is a natural map $\sigma :\Gamma 
_{r+1}\rightarrow \Gamma _{r}$, where $\sigma 
\left ( \delta _{i}\right )$ is defined by 
chopping off the two end edges of $\delta _{i}$. 
This map respects $\prec $ in the sense that, if 
$\delta \prec \delta '$, then $\sigma \left ( 
\delta \right )\prec \sigma \left ( \delta 
'\right )$ or $\sigma \left ( \delta \right 
)=\sigma \left ( \delta '\right )$. We can 
therefore choose the indexing so that, for every 
$j$, the set of those indices $i$ for which 
$\sigma \left ( \delta _{i}\right )=\gamma _{j}$ 
is of the form $k$, $k+1$, \dots , $k+l$. We will 
also denote by $\sigma $ the map $\left \lbrace 
1,\dots ,q\right \rbrace \rightarrow \left 
\lbrace 1,\dots ,p\right \rbrace $ defined by 
$\sigma \left ( \delta _{i}\right )=\gamma 
_{\sigma \left ( i\right )}$. 
\nobreak\par\nobreak 	For each $\delta _{i}$, let 
$h_{i}^{t}$ be the image under $\widetilde  
f_{t}'':\widetilde  S\rightarrow {\Bbb 
H}^{2}\subset {\Bbb H}^{3}$ of an 
$m_{t}$--geodesic of $\widetilde  S$ that is 
carried by $\widetilde  U$ and realizes $\delta 
_{i}$, if such a geodesic exists. Then,
$$A_{r+1}^{t}= R_{h_{1}^{t}}^{b_{t}\left ( \delta 
_{1}\right )} R_{h_{2}^{t}}^{b_{t}\left ( \delta 
_{2}\right )} \dots R_{h_{q}^{t}}^{b_{t}\left ( 
\delta _{q}\right )} \rho ''_{t}\left ( \xi 
\right ) .$$
Noting that $b_{t}\left ( \gamma _{j}\right 
)=\sum _{\sigma \left ( i\right )=j}b_{t}\left ( 
\delta _{i}\right )$, we can rewrite $A_{r}^{t}$ 
as
$$A_{r}^{t} =  R_{g_{\sigma (1)}^{t}}^{b_{t}\left 
( \delta _{1}\right )} R_{g_{\sigma 
(2)}^{t}}^{b_{t}\left ( \delta _{2}\right )} 
\dots R_{g_{\sigma (q)}^{t}}^{b_{t}\left ( \delta 
_{q}\right )} \rho ''_{t}\left ( \xi \right ) .$$
We conclude that
$$A_{r+1}^{t}-A_{r}^{t} = \sum _{i=1}^{q}  
R_{h_{1}^{t}}^{b_{t}\left ( \delta _{1}\right )} 
\dots R_{h_{i-1}^{t}}^{b_{t}\left ( \delta 
_{i-1}\right )} \left ( R_{h_{i}^{t}}^{b_{t}\left 
( \delta _{i}\right )} -R_{g_{\sigma 
(i)}^{t}}^{b_{t}\left ( \delta _{i}\right )} 
\right ) R_{g_{\sigma (i+1)}^{t}}^{b_{t}\left ( 
\delta _{i+1}\right )} \dots R_{g_{\sigma 
(q)}^{t}}^{b_{t}\left ( \delta _{q}\right )} \rho 
''_{t}\left ( \xi \right ) .$$
Similarly,
$$B_{r+1}^{t}-B_{r}^{t} = \sum _{i=1}^{q}  
R_{h_{1}^{t}}^{b_{t}'\left ( \delta _{1}\right )} 
\dots R_{h_{i-1}^{t}}^{b_{t}'\left ( \delta 
_{i-1}\right )} \left ( 
R_{h_{i}^{t}}^{b_{t}'\left ( \delta _{i}\right )} 
-R_{g_{\sigma (i)}^{t}}^{b_{t}'\left ( \delta 
_{i}\right )} \right ) R_{g_{\sigma 
(i+1)}^{t}}^{b_{t}'\left ( \delta _{i+1}\right )} 
\dots R_{g_{\sigma (q)}^{t}}^{b_{t}'\left ( 
\delta _{q}\right )} \rho ''_{t}\left ( \xi 
\right ) .$$
It follows that $C_{r+1}^{t}-C_{r}^{t}=\left ( 
A_{r+1}^{t}-A_{r}^{t}\right )-\left ( 
B_{r+1}^{t}-B_{r}^{t}\right )$ can be written as 
a sum of $q^{2}$ terms, each of the form
$$\multline
R_{h_{1}^{t}}^{b_{t}\left ( \delta _{1}\right )} 
\dots R_{h_{i-1}^{t}}^{b_{t}\left ( \delta 
_{i-1}\right )} \left ( R_{h_{i}^{t}}^{b_{t}\left 
( \delta _{i}\right )} -R_{g_{\sigma 
(i)}^{t}}^{b_{t}\left ( \delta _{i}\right )} 
\right )R_{g_{\sigma (i+1)}^{t}}^{b_{t}\left ( 
\delta _{i+1}\right )} \dots \\
 \dots R_{g_{\sigma (j-1)}^{t}}^{b_{t}\left ( 
\delta _{j-1}\right )} \left ( R_{g_{\sigma 
(j)}^{t}}^{b_{t}\left ( \delta _{j}\right )} 
-R_{g_{\sigma (j)}^{t}}^{b_{t}'\left ( \delta 
_{j}\right )} \right )R_{g_{\sigma 
(j+1)}^{t}}^{b_{t}'\left ( \delta _{j+1}\right )} 
\dots R_{g_{\sigma (q)}^{t}}^{b_{t}'\left ( 
\delta _{q}\right )} \rho ''_{t}\left ( \xi 
\right ) , 
\endmultline \tag 3$$

$$ 
R_{h_{1}^{t}}^{b_{t}\left ( \delta _{1}\right )} 
\dots R_{h_{i-1}^{t}}^{b_{t}\left ( \delta 
_{i-1}\right )} \left ( \left ( 
R_{h_{i}^{t}}^{b_{t}\left ( \delta _{i}\right )} 
-R_{g_{\sigma (i)}^{t}}^{b_{t}\left ( \delta 
_{i}\right )}\right )-\left ( 
R_{h_{i}^{t}}^{b_{t}'\left ( \delta _{i}\right )} 
-R_{g_{\sigma (i)}^{t}}^{b_{t}'\left ( \delta 
_{i}\right )} \right )\right )R_{g_{\sigma 
(i+1)}^{t}}^{b_{t}'\left ( \delta _{i+1}\right )} 
 \dots R_{g_{\sigma (q)}^{t}}^{b_{t}'\left ( 
\delta _{q}\right )} \rho ''_{t}\left ( \xi 
\right ) , 
 \tag 4$$
or
$$\multline
R_{h_{1}^{t}}^{b_{t}\left ( \delta _{1}\right )} 
\dots R_{h_{j-1}^{t}}^{b_{t}\left ( \delta 
_{j-1}\right )} \left ( R_{h_{j}^{t}}^{b_{t}\left 
( \delta _{j}\right )} 
-R_{h_{j}^{t}}^{b_{t}'\left ( \delta _{j}\right 
)} \right )R_{h_{j+1}^{t}}^{b_{t}'\left ( \delta 
_{j+1}\right )} \dots \\
 \dots R_{h_{i-1}^{t}}^{b_{t}'\left ( \delta 
_{i-1}\right )} \left ( 
R_{h_{i}^{t}}^{b_{t}'\left ( \delta _{i}\right )} 
-R_{g_{\sigma (i)}^{t}}^{b_{t}'\left ( \delta 
_{i}\right )} \right )R_{g_{\sigma 
(i+1)}^{t}}^{b_{t}'\left ( \delta _{i+1}\right )} 
\dots R_{g_{\sigma (q)}^{t}}^{b_{t}'\left ( 
\delta _{q}\right )} \rho ''_{t}\left ( \xi 
\right ) .
\endmultline \tag 5$$
\nobreak\par\nobreak 	To bound these terms, we 
will use the following estimate.

\proclaim {Lemma 8 }   Let $A_{1}$, $A_{2}$, 
\dots , $A_{n}$ be square matrices, and let the 
number $R$ bound the norm of all products 
$A_{i_{1}}A_{i_{2}}\dots A_{i_{p}}$ with 
$1\leqslant i_{1}<i_{2}<\dots <i_{p}\leqslant n$. 
Then, for every matrices $\varepsilon _{1}$, 
$\varepsilon _{2}$, \dots , $\varepsilon _{n}$, 
$$\left \Vert  \left ( A_{1}+\varepsilon 
_{1}\right ) \left ( A_{2}+\varepsilon _{2}\right 
) \dots \left ( A_{n}+\varepsilon _{n}\right 
)-A_{1}A_{2}\dots A_{n} \right \Vert  \leqslant  
R\left ( e^{nRE}-1\right )$$
where $E=\max_{i} \left \Vert \varepsilon 
_{i}\right \Vert $.  \endproclaim    
\demo {Proof}   If we expand $\left ( 
A_{1}+\varepsilon _{1}\right ) \left ( 
A_{2}+\varepsilon _{2}\right ) \dots \left ( 
A_{n}+\varepsilon _{n}\right )-A_{1}A_{2}\dots 
A_{n}$, each term in the expansion is the product 
of $k$ terms $\varepsilon _{i}$ and of $k+1$ 
terms $A_{i_{1}}A_{i_{2}}\dots A_{i_{s}}$ with 
$1\leqslant i_{1}<i_{2}<\dots <i_{s}\leqslant n$, 
for some $k$ between 1 and $n$. In addition, the 
number of terms with $k$ such $\varepsilon _{i}$ 
is equal to the binomial coefficient $\left ( {n 
\atop k}\right )$. It follows that
$$
\left \Vert  \left ( A_{1}+\varepsilon _{1}\right 
) \left ( A_{2}+\varepsilon _{2}\right ) \dots 
\left ( A_{n}+\varepsilon _{n}\right 
)-A_{1}A_{2}\dots A_{n} \right \Vert  \leqslant 
R\left ( \left ( 1+RE\right )^{n}-1\right )
								\leqslant  R\left ( e^{nRE}-1\right ).
$$
 \qed\enddemo 
 
\proclaim {Lemma 9 }   In the expressions 
\rom{(3)} to \rom{(5)}, the subterms of the form 
$R_{h_{k}^{t}}^{b_{t}\left ( \delta _{k}\right )} 
\dots R_{h_{l}^{t}}^{b_{t}\left ( \delta 
_{l}\right )} $, $R_{h_{k}^{t}}^{b_{t}'\left ( 
\delta _{k}\right )} \dots 
R_{h_{l}^{t}}^{b_{t}'\left ( \delta _{l}\right )} 
$, $R_{g_{\sigma (k)}^{t}}^{b_{t}\left ( \delta 
_{k}\right )} \dots R_{g_{\sigma 
(l)}^{t}}^{b_{t}\left ( \delta _{l}\right )}$, or 
$R_{g_{\sigma (k)}^{t}}^{b_{t}'\left ( \delta 
_{k}\right )} \dots R_{g_{\sigma 
(l)}^{t}}^{b_{t}'\left ( \delta _{l}\right )}$ 
are uniformly bounded (independent of $r$ and 
$t$).  \endproclaim    
\demo {Proof}   For every $i$ with $b_{t}\left ( 
\delta _{i}\right )\mathbin{\not =}0$, there is a 
leaf of the $m_{t}$--geodesic lamination 
underlying $b_{t}$ which realizes $\delta _{i}$, 
and we can consider its image $\bar h_{i}^{t}$ 
under $\widetilde  f_{t}''$. The main property we 
need is that the $\bar h_{i}^{t}$ are pairwise 
disjoint which, because the ordering of the 
$\delta _{i}$ is compatible with $\prec $, 
guarantees that $\bar h_{i}^{t}$ meets 
$\widetilde  f_{t}''\left ( \widetilde  c\right 
)$ closer to $\widetilde  f_{t}''\left ( 
\widetilde  x_{0}\right )$ than $\bar h_{i'}^{t}$ 
if $i<i'$. For $i_{1}<i_{2}<\dots <i_{p}$ with 
all $b_{t}\left ( \delta _{i_{j}}\right 
)\mathbin{\not =}0$, consider $R_{\bar 
h_{i_{1}}^{t}}^{b_{t}\left ( \delta 
_{i_{1}}\right )} \dots R_{\bar 
h_{i_{p}}^{t}}^{b_{t}\left ( \delta 
_{i_{p}}\right )} $. Because of the ordering of 
the intersections $\bar h_{i}^{t}\cap \widetilde  
f_{t}'\left ( \widetilde  c\right )$, the point 
$R_{\bar h_{i_{1}}^{t}}^{b_{t}\left ( \delta 
_{i_{1}}\right )} \dots R_{\bar 
h_{i_{p}}^{t}}^{b_{t}\left ( \delta 
_{i_{p}}\right )}  \widetilde  f_{t}'\left ( \xi 
\widetilde  x_{0}\right )$ can be connected to 
$\widetilde  f_{t}'\left ( \widetilde  
x_{0}\right )$ by a broken arc of the same length 
as $\widetilde  f_{t}'\left ( \widetilde  c\right 
)$. It follows that $R_{\bar 
h_{i_{1}}^{t}}^{b_{t}\left ( \delta 
_{i_{1}}\right )} \dots R_{\bar 
h_{i_{p}}^{t}}^{b_{t}\left ( \delta 
_{i_{p}}\right )} $ stays in a compact subset of 
the isometry group of ${\Bbb H}^{3}$; in 
particular, its norm is uniformly bounded by a 
constant $R>0$. 
\nobreak\par\nobreak 	Set $\varepsilon 
_{i}=R_{h_{i}^{t}}^{b_{t}\left ( \delta 
_{i}\right )} -R_{\bar h_{i}^{t}}^{b_{t}\left ( 
\delta _{i}\right )}$. Because 
$R_{h_{i}^{t}}^{b_{t}\left ( \delta _{i}\right 
)}$ and $R_{\bar h_{i}^{t}}^{b_{t}\left ( \delta 
_{i}\right )}$ are uniformly bounded, $\left 
\Vert \varepsilon _{i}\right \Vert $ is bounded 
by a constant times the distance between 
$h_{i}^{t}$ and $\bar h_{i}^{t}$. Because 
$h_{i}^{t}$ and $\bar h_{i}^{t}$ follow the same 
edge path of length $2r+1$, this distance is an 
$O\left ( e^{-Ar}\right )$ for some constant 
$A>0$ depending on $\widetilde  U$ and 
$\widetilde  c$. 
\nobreak\par\nobreak 	We are now in a position to 
apply Lemma~8. To prove that the product 
$R_{h_{k}^{t}}^{b_{t}\left ( \delta _{k}\right )} 
\dots R_{h_{l}^{t}}^{b_{t}\left ( \delta 
_{l}\right )} $ is uniformly bounded, Lemma~8 and 
the above estimate for $\varepsilon _{i}$ imply 
that it suffices to show that $\left ( l-k\right 
) e^{-Ar}$ is bounded. Although the number of 
edge paths $\delta \in \Gamma _{r+1}$ grows 
exponentially with $r$, the number of those for 
which $b_{t}\left ( \delta \right )\mathbin{\not 
=}0$ is bounded by a polynomial function of $r$ 
(this is a general fact about geodesic 
laminations, see for instance \cite{Bo2, 
Lemma~10}). It follows that $l-k=O\left ( 
r^{n}\right )$ for some $n$. As a consequence, 
$\left ( l-k\right ) e^{-Ar}$ is bounded. By 
Lemma~8, we conclude that all the products 
$R_{h_{k}^{t}}^{b_{t}\left ( \delta _{k}\right )} 
\dots R_{h_{l}^{t}}^{b_{t}\left ( \delta 
_{l}\right )} $ are uniformly bounded.
\nobreak\par\nobreak 	The proof of Lemma~9 for 
the products $R_{h_{k}^{t}}^{b_{t}'\left ( \delta 
_{k}\right )} \dots R_{h_{l}^{t}}^{b_{t}'\left ( 
\delta _{l}\right )} $, $R_{g_{\sigma 
(k)}^{t}}^{b_{t}\left ( \delta _{k}\right )} 
\dots R_{g_{\sigma (l)}^{t}}^{b_{t}\left ( \delta 
_{l}\right )}$ and $R_{g_{\sigma 
(k)}^{t}}^{b_{t}'\left ( \delta _{k}\right )} 
\dots R_{g_{\sigma (l)}^{t}}^{b_{t}'\left ( 
\delta _{l}\right )}$ is identical.  \qed\enddemo

\remark {Remark}   One could naively think that 
it is possible to greatly simplify the proof of 
Lemma~9 by taking $h_{i}^{t}=\bar h_{i}^{t}$ 
right away. However, it is not possible to do so 
simultaneously for the terms involving $b_{t}$ 
and those involving $b_{t}'$. In general, we 
cannot choose the $h_{i}^{t}$ so that $h_{i}^{t}$ 
is disjoint from $h_{i'}^{t}$ whenever 
$b_{t}\left ( \delta _{i}\right )b_{t}\left ( 
\delta _{i'}\right )\mathbin{\not =}0$ or 
$b_{t}'\left ( \delta _{i}\right )b_{t}'\left ( 
\delta _{i'}\right )\mathbin{\not =}0$. 
\endremark    

\medskip\nobreak 

\nobreak\par\nobreak 	We can now estimate 
$C_{r+1}^{t}-C_{r}^{t}$. 
\nobreak\par\nobreak 	In a term of type 
\rom{(3)}, the quantity 
$R_{h_{i}^{t}}^{b_{t}\left ( \delta _{i}\right )} 
-R_{g_{\sigma (i)}^{t}}^{b_{t}\left ( \delta 
_{i}\right )}$ is bounded by a constant times the 
distance from $h_{i}^{t}$ to $g_{\sigma \left ( 
i\right )}^{t}$, which is an $O\left ( 
e^{-Ar}\right )$ since these two geodesics follow 
the same edge path of length $2r+1$. The quantity 
$R_{g_{\sigma (j)}^{t}}^{b_{t}\left ( \delta 
_{j}\right )} -R_{g_{\sigma 
(j)}^{t}}^{b_{t}'\left ( \delta _{j}\right )}$ is 
bounded by a constant times $b_{t}\left ( \delta 
_{j}\right )-b_{t}'\left ( \delta _{j}\right )$. 
In \cite{Bo2, Lemma~2}, we give an explicit 
formula expressing $b_{t}\left ( \delta 
_{j}\right )$ in terms of the weights $b_{t}\left 
( e\right )$ it assigns to the edges $e$ of $U$. 
Because $\delta _{j}$ is an edge path of length 
$2r+3$, it follows from this formula that 
$b_{t}\left ( \delta _{j}\right )-b_{0}\left ( 
\delta _{j}\right )=O\left ( r\left \Vert 
b_{t}-b_{0}\right \Vert _{U}\right )=tO\bigl ( 
r\bigl \Vert \dot b_{0}\bigr \Vert _{U}\bigr ) $. 
Similarly,  $b_{t}\left ( \delta _{j}\right 
)=b_{0}\left ( \delta _{j}\right )+t\dot 
b_{0}\left ( \delta _{j}\right )=b_{0}\left ( 
\delta _{j}\right )+tO\bigl ( r\bigl \Vert \dot 
b_{0}\bigr \Vert _{U}\bigr ) $, and we conclude 
that $b_{t}\left ( \delta _{j}\right 
)-b_{t}'\left ( \delta _{j}\right )=tO\bigl ( 
r\bigl \Vert \dot b_{0}\bigr \Vert _{U}\bigr ) $. 
By Lemma~9, it follows that every term of type 
\rom{(3)} is of the form $tO\bigl ( re^{-Ar}\bigl 
\Vert \dot b_{0}\bigr \Vert _{U}\bigr ) $.
\nobreak\par\nobreak 	Similarly, every term of 
type \rom{(5)} is of the form $tO\bigl ( 
re^{-Ar}\bigl \Vert \dot b_{0}\bigr \Vert 
_{U}\bigr ) $.
\nobreak\par\nobreak 	In a term of type 
\rom{(4)}, the quantity $\bigl ( 
R_{h_{i}^{t}}^{b_{t}\left ( \delta _{i}\right )} 
-R_{g_{\sigma (i)}^{t}}^{b_{t}\left ( \delta 
_{i}\right )}\bigr ) -\bigl ( 
R_{h_{i}^{t}}^{b_{t}'\left ( \delta _{i}\right )} 
-R_{g_{\sigma (i)}^{t}}^{b_{t}'\left ( \delta 
_{i}\right )} \bigr ) $ is bounded by a constant 
times the product of $b_{t}\left ( \delta 
_{i}\right )-b_{t}'\left ( \delta _{i}\right )$ 
and of the distance from $h_{i}^{t}$ to 
$g_{\sigma \left ( i\right )}^{t}$. As above, we 
conclude that a term of type \rom{(4)} is of the 
form $tO\bigl ( re^{-Ar}\bigl \Vert \dot 
b_{0}\bigr \Vert _{U}\bigr ) $.
\nobreak\par\nobreak 	We saw that 
$C_{r+1}^{t}-C_{r}^{t}$ is a sum of $q^{2}$ terms 
of type \rom{(3)}, \rom{(4)} or \rom{(5)}. We 
also saw that $q=O\left ( r^{n}\right )$ for some 
$n$. Therefore, 
$$C_{r+1}^{t}-C_{r}^{t}=tO\bigl ( 
r^{2n+1}e^{-Ar}\bigl \Vert \dot b_{0}\bigr \Vert 
_{U}\bigr ) =tO\bigl ( e^{-A'r}\bigl \Vert \dot 
b_{0}\bigr \Vert _{U}\bigr ) $$
for any $A'>A$. This proves the first statement 
of Lemma~7.
\nobreak\par\nobreak 	The second statement of 
Lemma~7 is obtained by summing the differences 
$C_{r+1}^{t}-C_{r}^{t}$ from $r$ to $\infty $, 
since $\rho _{t}\left ( \xi \right )-\rho 
_{t}'\left ( \xi \right )=\lim_{r\rightarrow 
\infty } C_{r}^{t}$.  \qed\enddemo

\nobreak\par\nobreak 	Now, fix $r$ and let $t$ 
tend to $0^{+}$. For $i=1$, \dots , $p$,  $\bigl 
( R_{g_{i}^{t}}^{b_{t}\left ( \gamma _{i}\right 
)} -R_{g_{i}^{t}}^{b_{t}'\left ( \gamma 
_{i}\right )} \bigr ) /t=O\left ( b_{t}\left ( 
\gamma _{i}\right )-b_{t}'\left ( \gamma 
_{i}\right )\right )/t$. As $t$ tends to $0^{+}$, 
each $\left ( b_{t}\left ( \gamma _{i}\right 
)-b_{t}'\left ( \gamma _{i}\right )\right )/t$ 
converges to 0 since $b_{t}'\left ( \gamma 
_{i}\right )=b_{0}\left ( \gamma _{i}\right 
)+t\dot b_{0}\left ( \gamma _{i}\right )$. 
Therefore, for a fixed $r$, 
$$\align
C_{r}^{t} /t&= R_{g_{1}^{t}}^{b_{t}\left ( \gamma 
_{1}\right )} R_{g_{2}^{t}}^{b_{t}\left ( \gamma 
_{2}\right )} \dots R_{g_{p}^{t}}^{b_{t}\left ( 
\gamma _{p}\right )} \rho ''_{t}\left ( \xi 
\right )/t - R_{g_{1}^{t}}^{b_{t}'\left ( \gamma 
_{1}\right )} R_{g_{2}^{t}}^{b_{t}'\left ( \gamma 
_{2}\right )} \dots R_{g_{p}^{t}}^{b_{t}'\left ( 
\gamma _{p}\right )} \rho ''_{t}\left ( \xi 
\right )/t\\
&= \sum _{i=1}^{p} R_{g_{1}^{t}}^{b_{t}\left ( 
\gamma _{1}\right )} R_{g_{2}^{t}}^{b_{t}\left ( 
\gamma _{2}\right )} \dots 
R_{g_{i-1}^{t}}^{b_{t}\left ( \gamma _{i-1}\right 
)} {R_{g_{i}^{t}}^{b_{t}\left ( \gamma _{i}\right 
)} -R_{g_{i}^{t}}^{b_{t}'\left ( \gamma 
_{i}\right )} \over 
t}R_{g_{i+1}^{t}}^{b_{t}'\left ( \gamma 
_{i+1}\right )} \dots R_{g_{p}^{t}}^{b_{t}'\left 
( \gamma _{p}\right )} \rho ''_{t}\left ( \xi 
\right ) .
\endalign$$
 converges to 0 as $t$ tends to $0^{+}$.
\nobreak\par\nobreak 	It then follows from 
Lemma~7 that every limit point of $\left ( \rho 
_{t}\left ( \xi \right )-\rho _{t}'\left ( \xi 
\right )\right )/t$ as $t$ tends to $0^{+}$ is of 
the form $O\bigl ( e^{-Ar}\bigl \Vert \dot 
b_{0}\bigr \Vert _{U}\bigr ) $. 
\nobreak\par\nobreak 	This holds for every $r$. 
If we now let $r$ tend to $\infty $, we conclude 
that 0 is the only limit point of $\left ( \rho 
_{t}\left ( \xi \right )-\rho _{t}'\left ( \xi 
\right )\right )/t$ as $t$ tends to $0^{+}$, 
namely that the two curves $t\mapsto \rho 
_{t}\left ( \xi \right )$ and $t\mapsto \rho 
_{t}'\left ( \xi \right )\in \roman 
{Isom}^{+}\left ( {\Bbb H}^{3}\right )$ are 
tangent at $t=0$.
\nobreak\par\nobreak 	This shows that the two 
curves $t\mapsto \rho _{t}$ and $t\mapsto \rho 
_{t}'\in {\Cal R}\left ( S\right )$ are tangent 
at $t=0$. As a consequence, $t\mapsto \rho _{t}$ 
has a tangent vector $\dot \rho _{0}$ at $t=0$, 
which is equal to $\dot \rho _{0}'=T_{\left ( 
m_{0}, \bar b_{0}\right )} \varphi _{\lambda 
}\bigl ( \dot m_{0},\dot b_{0}\bigr ) $.
\nobreak\par\nobreak 	This concludes the proof of 
Proposition~5 under the additional assumption 
that the geodesic laminations underlying the 
$b_{t}$ converge to some sublamination of 
$\lambda $.
\nobreak\par\nobreak 	In the general case, let 
$t_{n}$, $n\in {\Bbb N}$, be a sequence 
converging to $0^{+}$, such that the geodesic 
lamination underlying $b_{t_{n}}$ converges to 
some lamination $\lambda '$ for the Hausdorff 
topology. The geodesic lamination $\lambda '$ 
must contain the supports of $b_{0}$ and $\dot 
b_{0}$. We can therefore consider 
$b_{t}'=b_{0}+t\dot b_{0}$ as a transverse 
cocycle for $\lambda '$ as well as for $\lambda 
$; the same holds for its reduction $\bar b_{t}'$ 
modulo $2\pi $. Note that $\varphi _{\lambda 
'}\bigl ( m_{t}, \bar b_{t}'\bigr ) =\varphi 
_{\lambda }\bigl ( m_{t}, \bar b_{t}'\bigr ) 
=\rho _{t}'$. Then, the same argument as above 
shows that the ``discrete curve'' $t_{n}\mapsto 
\rho _{t_{n}}$ is tangent to the curve $t\mapsto 
\rho _{t}'$ at 0, in the sense that 
$\lim_{n\rightarrow \infty } \left ( \rho 
_{t_{n}}\left ( \xi \right )-\rho _{t_{n}}'\left 
( \xi \right )\right )/t_{n}=0$ for every $\xi 
\in \pi _{1}\left ( S\right )$.  Since this 
property holds for any such subsequence $t_{n}$, 
$n\in {\Bbb N}$, this shows that the two curves 
$t\mapsto \rho _{t}$ and $t\mapsto \rho _{t}'$ 
are tangent at $t=0$. Again, it follows that 
$t\mapsto \rho _{t}$ has a tangent vector $\dot 
\rho _{0}$ at $t=0$ which is equal to $\dot \rho 
_{0}'=T_{\left ( m_{0}, \bar b_{0}\right )} 
\varphi _{\lambda }\bigl ( \dot m_{0},\dot 
b_{0}\bigr ) $, and this completes the proof of 
Proposition~5.
\qed\enddemo

\nobreak\par\nobreak 	By Proposition~5, the map 
$\varphi :{\Cal T}\left ( S\right )\times {\Cal 
M}{\Cal L}\left ( S\right )\rightarrow {\Cal 
R}\left ( S\right )$  has a tangent map $T_{\left 
( m,b\right )}\varphi :T_{m}{\Cal T}\left ( 
S\right )\times T_{b}{\Cal M}{\Cal L}\left ( 
S\right )\rightarrow T_{\varphi \left ( m,b\right 
)}{\Cal R}\left ( S\right )$ everywhere. If, in 
addition, the support of $b$ is a maximal 
geodesic lamination $\lambda $, then $T_{b}{\Cal 
M}{\Cal L}\left ( S\right )\cong {\Cal H}\left ( 
\lambda ;{\Bbb R}\right )$ and $T_{\left ( 
m,b\right )}\varphi =T_{\left ( m,b\right 
)}\varphi _{\lambda }$. Since $\varphi _{\lambda 
}$ is a local diffeomorphism, this immediately 
shows that $T_{\left ( m,b\right )}\varphi $ is 
invertible when the support of $b$ is a maximal 
geodesic lamination. The general case requires 
more work.

\head {\S }3. Proof that $T_{\left ( m,b\right 
)}\varphi :T_{m}{\Cal T}\left ( S\right )\times 
T_{b}{\Cal M}{\Cal L}\left ( S\right )\rightarrow 
T_{\varphi \left ( m,b\right )}{\Cal R}\left ( 
S\right )$ is injective\endhead    

\proclaim {Proposition 10}   The tangent map 
$T_{\left ( m,b\right )}\varphi :T_{m}{\Cal 
T}\left ( S\right )\times T_{b}{\Cal M}{\Cal 
L}\left ( S\right )\rightarrow T_{\varphi \left ( 
m,b\right )}{\Cal R}\left ( S\right )$ is 
injective. \endproclaim    
\demo {Proof}   Let $v'=\bigl ( \dot m', \dot 
b'\bigr ) $ and $v''=\bigl ( \dot m'', \dot 
b''\bigr ) $ be two tangent vectors at $\left ( 
m_{0}, b_{0}\right )$ such that $T_{\left ( 
m_{0},b_{0}\right )}\varphi \left ( v'\right 
)=T_{\left ( m_{0},b_{0}\right )}\varphi \left ( 
v''\right )$. We want to show that $v'=v''$.
\nobreak\par\nobreak 	By Proposition~5, $T_{\left 
( m_{0},b_{0}\right )}\varphi \left ( v'\right 
)=T_{\left ( m_{0}, b_{0}\right )}\varphi 
_{\lambda '}\bigl ( \dot m', \dot b'\bigr ) $ 
where $\lambda '$ is any maximal geodesic 
lamination containing the supports of $b_{0}$ and 
$\dot b$. Similarly, $T_{\left ( m,b\right 
)}\varphi \left ( v''\right )=T_{\left ( m_{0}, 
b_{0}\right )}\varphi _{\lambda ''}\bigl ( \dot 
m'', \dot b''\bigr ) $ where $\lambda ''$ is any 
maximal geodesic lamination containing the 
supports of $b_{0}$ and $\dot b''$. 
\proclaim {Lemma 11}   The support of $\dot b'$ 
does not cross the support of $\dot b''$.  
\endproclaim    
\demo {Proof}   Suppose that there is a leaf $g'$ 
of the support of $\dot b'$ that intersects 
transversely in $x$ a leaf $g''$ of the support 
of $\dot b''$. Without loss of generality, we may 
assume that $g'$ is in the boundary of $S-\lambda 
'$ and that $g''$ is in the boundary of 
$S-\lambda ''$. 
\nobreak\par\nobreak 	Let $\rho _{t}\in {\Cal 
R}\left ( S\right )$, $t\in \left [0,\varepsilon 
\right [$, be a family of representations with 
$\rho _{0}=\varphi _{\lambda '}\left ( 
m_{0},b_{0}\right )=\varphi _{\lambda ''}\left ( 
m_{0},b_{0}\right )$ and $\dot \rho _{0}= 
T_{\left ( m_{0},b_{0}\right )}\varphi _{\lambda 
'}\bigl ( \dot m', \dot b'\bigr ) = T_{\left ( 
m_{0},b_{0}\right )}\varphi _{\lambda ''}\bigl ( 
\dot m'', \dot b''\bigr ) $. For $t$ small 
enough, $\rho _{t}$ determines a pleated surface 
$f_{t}'=\bigl ( \widetilde  f_{t}',\rho _{t}\bigr 
) $ with pleating locus $\lambda '$ and a pleated 
surface $f_{t}''=\bigl ( \widetilde  f_{t}'',\rho 
_{t}\bigr ) $ with pleating locus $\lambda ''$. 
Let $m_{t}'\in {\Cal T}\left ( S\right )$ and 
$b_{t}'\in {\Cal H}\left ( \lambda ';{\Bbb 
R}/2\pi {\Bbb Z}\right )$ (resp. $m_{t}''\in 
{\Cal T}\left ( S\right )$ and $b_{t}''\in {\Cal 
H}\left ( \lambda '';{\Bbb R}/2\pi {\Bbb Z}\right 
)$) be the pull back metric and the bending 
cocycles of $f_{t}'$ (resp. $f_{t}''$). Namely, 
$\rho _{t}=\varphi _{\lambda '}\left ( 
m_{t}',b_{t}'\right )=\varphi _{\lambda ''}\left 
( m_{t}'',b_{t}''\right )$. Note that 
$b_{0}'=b_{0}''=b_{0}$, $\dot b_{0}'=\dot b'$ and 
$\dot b_{0}''=\dot b''$.
\nobreak\par\nobreak 	Lift $x$ to a point 
$\widetilde  x$ in the universal covering 
$\widetilde  S$, and let $\widetilde  g'$ and 
$\widetilde  g''$ be the lifts of $g'$ and $g''$ 
passing through $\widetilde  x$, respectively. We 
want to compare the respective positions of the 
geodesics $\widetilde  f_{t}'\left ( \widetilde  
g_{t}'\right )$ and $\widetilde  f_{t}''\left ( 
\widetilde  g_{t}''\right )$ of ${\Bbb H}^{3}$, 
where $\widetilde  g_{t}'$ is the 
$m_{t}'$--geodesic of $\widetilde  S$ 
corresponding to $\widetilde  g'$ and $\widetilde 
 g_{t}''$ is the $m_{t}''$--geodesic 
corresponding to $\widetilde  g''$. Because 
$\widetilde  f_{0}'=\widetilde  f_{0}''$, the 
geodesics $\widetilde  f_{0}'\left ( \widetilde  
g_{0}'\right )$ and $\widetilde  f_{0}''\left ( 
\widetilde  g_{0}''\right )$ are coplanar and 
meet in one point. 
\nobreak\par\nobreak 	 If $\widehat g_{t}''$ 
denotes the $m_{t}'$--geodesic corresponding to 
$\widetilde  g''$, $\widetilde  f_{t}''\left ( 
\widetilde  g_{t}''\right )$ is also the geodesic 
of ${\Bbb H}^{3}$ that is asymptotic to 
$\widetilde  f_{t}'\left ( \widehat g_{t}''\right 
)$. Because $g'$ and $g''$ intersect, they have 
to be disjoint from the support of $b_{0}$. This 
implies that $\dot b'\left ( k''\right )\geqslant 
0$ for every arc $k''$ contained in $g''$. 
Indeed, $\dot b'\in T_{b_{0}}{\Cal M}{\Cal 
L}\left ( S\right )$ is tangent to a family of 
measured laminations $b_{t}\in {\Cal M}{\Cal 
L}\left ( S\right )$ with $b_{0}\left ( k\right 
)=0$ and $b_{t}\left ( k\right )\geqslant 0$; 
compare \cite{Bo2, Theorem~19}. It follows that, 
infinitesimally, $\widetilde  f_{t}'\left ( 
\widehat g_{t}''\right )$ bends everywhere in the 
direction of the negative side of $\widetilde  
f_{0}'\bigl ( \widetilde  S\bigr ) $. 
Intuitively, this will imply that, as $t$ moves 
away from 0, $\widetilde  f_{t}''\left ( 
\widetilde  g_{t}''\right )$ moves away from 
$\widetilde  f_{t}'\left ( \widetilde  
g_{t}'\right )$ in the direction of the negative 
side of $\widetilde  f_{0}'\bigl ( \widetilde  
S\bigr ) $. We need to quantify this.
\nobreak\par\nobreak 	By \cite{Bo3, 
Corollary~32}, for every component $P$ of 
$\widetilde  S-\widetilde  \lambda '$, the 
infinite triangle $\widetilde  f_{t}'\left ( 
P\right )\subset {\Bbb H}^{3}$ depends 
differentiably on the representation $\rho _{t}$. 
By our assumption that $g$ is a boundary leaf, it 
follows that $\widetilde  f_{t}'\left ( 
\widetilde  g_{t}'\right )$ depends 
differentiably on $\rho _{t}$. Since the same 
property holds for $\widetilde  f_{t}''\left ( 
\widetilde  g_{t}''\right )$, the length $l_{t}$ 
of the shortest geodesic arc from $\widetilde  
f_{t}'\left ( \widetilde  g_{t}'\right )$ to 
$\widetilde  f_{t}''\left ( \widetilde  
g_{t}''\right )$ also depends differentiably on 
$\rho _{t}$.
\nobreak\par\nobreak 	To estimate the derivative 
$\dot l_{0}$, normalize $\rho _{t}$ and 
$\widetilde  f_{t}'$ so that $\widetilde  f_{t}'$ 
sends the component of $\widetilde  S-\widetilde  
\lambda $ that is adjacent to $\widetilde  g'$ to 
a fixed ideal triangle in ${\Bbb H}^{2}\subset 
{\Bbb H}^{3}$. Then, $\widetilde  f_{t}'\left ( 
\widehat g_{t}''\right )$ is obtained from the 
geodesic $\widetilde  f_{0}'\left ( \widetilde  
g_{0}''\right )\subset {\Bbb H}^{2}$ by, first 
moving it in ${\Bbb H}^{2}$ to reflect the 
passage from the metric $m_{0}'$ to $m_{t}'$, and 
then bending this geodesic by successive 
rotations along geodesics of ${\Bbb H}^{2}$, 
following a formula analogous to \rom{(1)}. Let 
$\widetilde  h''$ be a half-line in $\widetilde  
g''$ which crosses the support of $\dot b'$ and 
which originates in the component of $\widetilde  
S-\widetilde  \lambda $ that is adjacent to 
$\widetilde  g'$; we will denote by $\widetilde  
h_{t}''$, $\widehat h_{t}''$ the subsets of 
$\widetilde  g_{t}''$, $\widehat g_{t}''$ 
corresponding to $\widetilde  h''$. Let $\theta 
_{t}^{+}$ be the visual amount by which the end 
point of $\widetilde  f_{t}'\bigl ( \widehat 
h_{t}''\bigr ) $ dips below ${\Bbb H}^{2}$, as 
measured from a fixed base point on ${\Bbb 
H}^{2}$. 
\nobreak\par\nobreak 	The derivative of $\theta 
_{t}^{+}$ at $t=0$ is given by the formula
$$\dot \theta _{0}^{+} = \displaystyle \int 
_{\widetilde  f_{0}'\left ( \widetilde  
h_{0}''\right )} A^{+}\left ( u\right )\thinspace 
d\dot b'\left ( u\right ) $$
where:  $d\dot b'$ is the distribution induced by 
$\dot b'$ on $\widetilde  f_{0}'\bigl ( 
\widetilde  h_{0}''\bigr ) $, which is actually a 
(countably additive) measure since $\dot b'\left 
( k''\right )\geqslant 0$ for every arc $k''$ 
contained in $g''$; $A^{+}\left ( u\right )>0$ 
denotes the amount by which the end point of 
$\widetilde  f_{0}'\bigl ( \widetilde  
h_{0}''\bigr ) $ dips under ${\Bbb H}^{2}$ when 
we apply to it the infinitesimal rotation around 
the leaf of $\widetilde  f_{0}'\bigl ( \widetilde 
 \lambda '\bigr ) $ passing through $u\in 
\widetilde  f_{0}'\bigl ( \widetilde  
h_{0}''\bigr ) $, if it exists. This formula is 
easily obtained by formal computations. To 
justify these formal computations (and show that 
the integral really converges), it suffices to 
note that $-\log A^{+}\left ( u\right )$ is at 
least a constant times the distance from $u$ to 
the base point and that, for every arc $k$ of 
length $\geqslant 1$ in $\widetilde  f_{0}'\bigl 
( \widetilde  h_{0}''\bigr ) $, $\dot b\left ( 
k\right )$ is bounded by a constant (depending on 
$\dot b$ but not on $k$) times the length of $k$. 
\nobreak\par\nobreak 	The important part here is 
that $\dot \theta _{0}^{+}>0$, which holds 
because $\widetilde  h''$ crosses the support of 
$\dot b$. A similar formula gives that $\dot 
\theta _{0}^{-}\geqslant 0$, where $\theta 
_{t}^{-}$ denotes the visual amount by which the 
other end point of $\widetilde  f_{t}'\left ( 
\widehat g_{t}''\right )$ dips below ${\Bbb 
H}^{2}$. Combining these two properties, it 
follows that $\dot l_{0}>0$.
\nobreak\par\nobreak 	This proves that, for 
$t>0$, the shortest geodesic arc from $\widetilde 
 f_{t}'\left ( \widetilde  g_{t}'\right )$ to 
$\widetilde  f_{t}''\left ( \widetilde  
g_{t}''\right )$ is non-trivial and points in the 
direction of the negative side of $\widetilde  
f_{0}'\bigl ( \widetilde  S\bigr ) $. But the 
argument is symmetric. Exchanging primes and 
double primes, we obtain that, for $t>0$, the 
opposite shortest geodesic arc from $\widetilde  
f_{t}''\left ( \widetilde  g_{t}''\right )$ to 
$\widetilde  f_{t}'\left ( \widetilde  
g_{t}'\right )$ must also point in the direction 
of the negative side of $\widetilde  f_{0}''\bigl 
( \widetilde  S\bigr ) =\widetilde  f_{0}'\bigl ( 
\widetilde  S\bigr ) $, a contradiction. 
\qed\enddemo

\nobreak\par\nobreak 	By Lemma~11, the supports 
of $\dot b'$ and $\dot b''$ do not cross each 
other. Therefore, there exists a maximal geodesic 
lamination $\lambda $ which contains the supports 
of $b_{0}$, $\dot b'$ and $\dot b''$. As a 
consequence, we can choose our geodesic 
laminations $\lambda '$, $\lambda ''$ so that 
$\lambda '=\lambda ''=\lambda $. 
\nobreak\par\nobreak 	Then,
$$T_{\left ( m_{0}, b_{0}\right )}\varphi 
_{\lambda }\left ( v'\right )=T_{\left ( m_{0}, 
b_{0}\right )}\varphi \left ( v'\right )=T_{\left 
( m_{0}, b_{0}\right )}\varphi \left ( v''\right 
)=T_{\left ( m_{0}, b_{0}\right )}\varphi 
_{\lambda }\left ( v''\right ).$$
Since $\varphi _{\lambda }$ is a diffeomorphism, 
its tangent map is a linear isomorphism, and it 
follows that $v'=v''$.  \qed\enddemo

\head {\S }4. Proof of Theorems~2 and 3\endhead   
\nobreak\par\nobreak 	Theorems~2 and 3 
immediately follow from Lemma~4, Corollary 6 and 
Proposition~10.
\nobreak\par\nobreak 	Indeed, for a connected 
surface $S$ of finite type and negative Euler 
characteristic,  the map $\varphi :{\Cal T}\left 
( S\right )\times {\Cal M}{\Cal L}\left ( S\right 
)\rightarrow {\Cal R}\left ( S\right )$ is the 
composition of the Thurston homeomorphism $\psi 
:{\Cal T}\left ( S\right )\times {\Cal M}{\Cal 
L}\left ( S\right )\rightarrow {\Cal P}\left ( 
S\right )$ and of the monodromy map $\theta 
:{\Cal P}\left ( S\right )\rightarrow {\Cal 
R}\left ( S\right )$. Because $\theta $ is a 
local diffeomorphism, $\varphi $ is a local 
homeomorphism. By Corollary~6, $\varphi $ admits 
a tangent map everywhere, and Proposition~10 
shows that this tangent map is injective. From 
Lemma~4, we conclude that any local inverse 
$\varphi ^{-1}$ for $\varphi $ is also 
tangentiable. Because $\theta $ is a local 
diffeomorphism, this shows that $\psi $ and $\psi 
^{-1}$ are tangentiable. This proves Theorem~3.
\nobreak\par\nobreak 	For a hyperbolic 
3--manifold $M$, the map $\mu \times \beta :{\Cal 
Q}{\Cal D}\left ( M\right )\rightarrow {\Cal 
T}\left ( \partial C_{M}\right )\times {\Cal 
M}{\Cal L}\left ( \partial C_{M}\right )$ locally 
coincides with the composition $\varphi 
^{-1}\circ R$ near the metric $M$ where, as in 
the introduction, ${\Cal R}\left ( \partial 
C_{M}\right )$ denotes the product $\prod 
_{i=1}^{n}{\Cal R}\left ( S_{i}\right )$ of the 
representation spaces corresponding to the 
components $S_{1}$, \dots , $S_{n}$ of $\partial 
C_{M}$, where $R:{\Cal Q}{\Cal D}\left ( M\right 
)\rightarrow {\Cal R}\left ( \partial C_{M}\right 
)$ is defined by restriction of the holonomy map, 
where $\varphi :{\Cal T}\left ( \partial 
C_{M}\right )\times {\Cal M}{\Cal L}\left ( 
\partial C_{M}\right )\rightarrow {\Cal R}\left ( 
\partial C_{M}\right )$ is defined as the product 
of the bending maps $\varphi _{i}:{\Cal T}\left ( 
S_{i}\right )\times {\Cal M}{\Cal L}\left ( 
S_{i}\right )\rightarrow {\Cal R}\left ( 
S_{i}\right )$, and where $\varphi ^{-1}$ is the 
local inverse defined near the representation 
$R\left ( M\right )$ and $\left ( \mu \left ( 
M\right ),\beta \left ( M\right )\right )$. As 
above, a combination of Corollary~6, 
Proposition~10 and Lemma~4 shows that each local 
inverse $\varphi _{i}^{-1}$ is tangentiable. 
Therefore, the local inverse $\varphi ^{-1}$ is 
tangentiable. Since $R$ is a differentiable map 
between differentiable manifolds, it follows that 
$\mu \times \beta $ is tangentiable. Composing 
with the (clearly tangentiable) projection 
$P:{\Cal T}\left ( \partial C_{M}\right )\times 
{\Cal M}{\Cal L}\left ( \partial C_{M}\right 
)\rightarrow {\Cal M}{\Cal L}\left ( \partial 
C_{M}\right )$, we conclude that $\beta $ is 
tangentiable everywhere. This proves Theorem~2.
\nobreak\par\nobreak 	The same argument shows 
that $\mu $ is tangentiable everywhere. To show 
that $\mu $ is continuously differentiable in the 
usual sense, we have to show that its tangent 
maps are linear and vary continuously with their 
base point. This will be done in the next 
section.

\head {\S }5. Proof of Theorem~1\endhead    
\nobreak\par\nobreak 	By the same arguments as in 
{\S }4, Theorem~1 immediately follows from the 
folllowing result.

\proclaim {Proposition 12  }   Let $S$ be a 
connected oriented surface of finite type and 
negative Euler characteristic. Then the 
composition $Q\circ \varphi ^{-1}$ of any local 
inverse $\varphi ^{-1}$ for the bending map 
$\varphi :{\Cal T}\left ( S\right )\times {\Cal 
M}{\Cal L}\left ( S\right )\rightarrow {\Cal 
R}\left ( S\right )$ and of the projection 
$Q:{\Cal T}\left ( S\right )\times {\Cal M}{\Cal 
L}\left ( S\right )\rightarrow {\Cal T}\left ( 
S\right )$ is continuously differentiable. 
\endproclaim    
\demo {Proof}   Let $\rho _{0}\in {\Cal R}\left ( 
S\right )$ and $\left ( m_{0},b_{0}\right 
)=\varphi ^{-1}\left ( \rho _{0}\right )$. By 
Corollary~6, Proposition~10 and Lemma~4, $\varphi 
^{-1}$ has a tangent map at $\rho _{0}$ and 
$T_{\rho _{0}}\varphi ^{-1}=\left ( T_{\left ( 
m_{0},b_{0}\right )}\varphi \right )^{-1}$. 
\nobreak\par\nobreak 	By Corollary~6, the 
restriction of $T_{\left ( m_{0},b_{0}\right 
)}\varphi $ to $T_{m_{0}}{\Cal T}\left ( S\right 
)\times 0$ coincides with the restriction of 
$T_{\left ( m_{0},b_{0}\right )}\varphi _{\lambda 
}$ for any maximal geodesic lamination $\lambda $ 
containing the support of $b$. In particular, 
this restriction of $T_{\left ( m_{0},b_{0}\right 
)}\varphi $ to $T_{m_{0}}{\Cal T}\left ( S\right 
)\times 0$ is linear. Let $P_{\rho _{0}}\subset 
T_{\rho _{0}}{\Cal R}\left ( S\right )$ denote 
the linear subspace $T_{\left ( m_{0},b_{0}\right 
)}\varphi \left ( T_{m_{0}}{\Cal T}\left ( 
S\right )\times 0\right )$; note that $P_{\rho 
_{0}}$ depends on $\rho _{0}$, but also on the 
choice of the local inverse $\varphi ^{-1}$.
\nobreak\par\nobreak 	To consider the image of 
$0\times T_{b_{0}}{\Cal M}{\Cal L}\left ( S\right 
)$ under $T_{\left ( m_{0},b_{0}\right )}\varphi 
$, we will exploit the complex structure of 
${\Cal R}\left ( S\right )$ coming from the 
complex structure of the group $\roman 
{Isom}^{+}\left ( {\Bbb H}^{3}\right )=\roman 
{PSL}_{2}\left ( {\Bbb C}\right )$. Indeed, it is 
showed in \cite{Bo3, {\S }10} that, for every 
maximal geodesic lamination $\lambda $ containing 
the support of $b$, the differential $T_{\left ( 
m_{0},b_{0}\right )}\varphi _{\lambda }$ sends 
$0\times {\Cal H}\left ( \lambda ;{\Bbb R}\right 
)$ to the subspace $iP_{\rho _{0}}$ obtained from 
$P_{\rho _{0}}$ by multiplication by $i$; see 
also the proof of Lemma~13 below. By Corollary~6, 
this implies that $T_{\left ( m_{0},b_{0}\right 
)}\varphi $ sends $0\times T_{b_{0}}{\Cal M}{\Cal 
L}\left ( S\right )$ inside $iP_{\rho _{0}}$. 
Because $T_{\left ( m_{0},b_{0}\right )}\varphi $ 
is invertible, the image of $0\times 
T_{b_{0}}{\Cal M}{\Cal L}\left ( S\right )$ by 
$T_{\left ( m_{0},b_{0}\right )}\varphi $ is 
actually equal to $iP_{\rho _{0}}$. (As an aside, 
since $T_{\left ( m_{0},b_{0}\right )}\varphi $ 
identifies $0\times T_{b_{0}}{\Cal M}{\Cal 
L}\left ( S\right )$ to $iP_{\rho _{0}}$, this 
defines on $T_{b_{0}}{\Cal M}{\Cal L}\left ( 
S\right )$ a linear structure which is compatible 
with the linear structures of the faces and 
depends only on $m_{0}$).
\nobreak\par\nobreak 	We can then compute the 
tangent map $T_{\rho _{0}}\left ( Q\circ \varphi 
^{-1}\right ):T_{\rho _{0}}{\Cal R}\left ( 
S\right )\rightarrow T_{m_{0}}{\Cal T}\left ( 
S\right )$. By Corollary~6, $T_{\rho _{0}}\left ( 
Q\circ \varphi ^{-1}\right )$ is just the 
composition $\Phi _{\rho _{0}}^{-1}\circ \Pi 
_{\rho _{0}}$ of the projection $\Pi _{\rho 
_{0}}$ of $T_{\rho _{0}}{\Cal R}\left ( S\right 
)$ onto $P_{\rho _{0}}$ parallel to $iP_{\rho 
_{0}}$ and of the inverse of the linear 
isomorphism $\Phi _{\rho _{0}}:T_{m_{0}}{\Cal 
T}\left ( S\right )\rightarrow P_{\rho _{0}}$ 
induced by $T_{\left ( m_{0},b_{0}\right 
)}\varphi $. In particular, $T_{\rho _{0}}\left ( 
Q\circ \varphi ^{-1}\right )$ is linear, and 
$Q\circ \varphi ^{-1}$ is differentiable in the 
usual sense.
\nobreak\par\nobreak 	It remains to show that 
$T_{\rho _{0}}\left ( Q\circ \varphi ^{-1}\right 
)$ depends continuously on $\rho _{0}$. 

\proclaim {Lemma 13 }   The linear map $\Phi 
_{\rho _{0}}:T_{m_{0}}{\Cal T}\left ( S\right 
)\rightarrow P_{\rho _{0}}$ depends continuously 
on $\rho _{0}$. \endproclaim    
\demo {Proof}   We will again make use of the 
complex structure of ${\Cal R}\left ( S\right )$. 
\nobreak\par\nobreak 	If $\lambda $ is a maximal 
geodesic lamination containing the support of 
$b_{0}$, we saw that $\varphi _{\lambda }$ 
provides a local parametrization of ${\Cal 
R}\left ( S\right )$ near $\rho _{0}$. This 
parametrization associates to each representation 
near $\rho _{0}$ the pull back metric $m_{\rho 
}\in {\Cal T}\left ( S\right )$ and the bending 
cocycle $b_{\rho }\in {\Cal H}\left ( \lambda 
;{\Bbb R}/2\pi {\Bbb Z}\right )$ of the pleated 
surface with pleating locus $\lambda $ 
corresponding to $\rho $. In \cite{Bo3}, we also 
associated to $m_{\rho }$ on $S$ a {\bi shearing 
cocycle \/} $s_{\rho }\in {\Cal H}\left ( \lambda 
;{\Bbb R}\right )$, and we combined $s_{\rho }$ 
and $b_{\rho }$ into a complex cocycle $s_{\rho 
}+ib_{\rho }\in {\Cal H}\left ( \lambda ;{\Bbb 
C}/2\pi i{\Bbb Z}\right )$. We showed that this 
provides a biholomorphic parametrization of a 
neighborhood of $\rho _{0}$ by an open subset of 
${\Cal H}\left ( \lambda ;{\Bbb C}/2\pi i{\Bbb 
Z}\right )$.
\nobreak\par\nobreak 	If $U$ is a train track 
carrying $\lambda $, each transverse cocycle 
$a\in {\Cal H}\left ( \lambda ;{\Bbb C}/2\pi 
i{\Bbb Z}\right )$ associates to each edge $e$ of 
$U$ a weight $a\left ( e\right )\in {\Bbb C}/2\pi 
i{\Bbb Z}$.  This defines a linear isomorphism 
between ${\Cal H}\left ( \lambda ;{\Bbb C}/2\pi 
i{\Bbb Z}\right )$ and the space ${\Cal W}\left ( 
U;{\Bbb C}/2\pi i{\Bbb Z}\right )$ of all such 
systems of edge weights that satisfy the 
classical {\bi switch relations \/}, namely such 
that, at each switch of $U$, the sum of the 
weights of the edges coming on one side is equal 
to the sum of the weights of the edges coming on 
the other side; see for instance \cite{Bo1}. 
\nobreak\par\nobreak 	Combining these two 
parametrizations, we get a holomorphic map $\psi 
_{\lambda }:{\Cal U}\rightarrow {\Cal R}\left ( 
S\right )$ which restricts to a homeomorphism 
between an open subset ${\Cal U}$ of ${\Cal 
W}\left ( U;{\Bbb C}/2\pi i{\Bbb Z}\right )$ and 
a neighborhood $\psi _{\lambda }\left ( {\Cal 
U}\right )$ of $\rho _{0}$. 
\nobreak\par\nobreak 	The main point of using 
edge weights instead of transverse cocycles is 
that we can compare these maps as we vary the 
geodesic lamination $\lambda $. If $\lambda 
_{n}$, $n\in {\Bbb N}$, is a sequence of geodesic 
lamination that converges to $\lambda $ for the 
Hausdorff topology as $n$ tends to $\infty $, the 
estimates of \cite{Bo3, {\S }4} show that, for 
$n$ large enough, the $\psi _{\lambda _{n}}$ are 
also defined on the same ${\Cal U}\subset {\Cal 
W}\left ( U;{\Bbb C}/2\pi i{\Bbb Z}\right )$ and 
uniformly converge to $\psi _{\lambda }$ on 
${\Cal U}$. Because the $\psi _{\lambda _{n}}$ 
are holomorphic, we also have uniform convergence 
of their tangent maps. We conclude that if, in 
addition, we have a sequence of edge weight 
systems $A_{n}\in {\Cal U}$ converging to some 
$A\in {\Cal U}$ and a sequence of tangent vectors 
$\dot A_{n}\in T_{A_{n}}{\Cal U}={\Cal W}\left ( 
U;{\Bbb C}\right )$ converging to $\dot A\in 
T_{A}{\Cal U}={\Cal W}\left ( U;{\Bbb C}\right )$ 
then, in ${\Cal R}\left ( S\right )$, the tangent 
vectors $T_{A_{n}}\psi _{\lambda _{n}}\bigl ( 
\dot A_{n}\bigr ) $ converge to $T_{A}\psi 
_{\lambda }\bigl ( \dot A\bigr ) $ as $n$ tends 
to $\infty $. 
\nobreak\par\nobreak 	If we restrict attention to 
real cocycles (and consequently to totally 
geodesic pleated surfaces and Fuchsian 
representations), we similarly have a real 
analytic map $\theta _{\lambda }:{\Cal 
V}\rightarrow {\Cal T}\left ( S\right )$ which 
restricts to a homeomoprhism between an open 
subset ${\Cal V}$ of ${\Cal W}\left ( U;{\Bbb 
R}\right )$ and a neighborhood $\theta _{\lambda 
}\left ( {\Cal V}\right )$ of $m_{0}\in {\Cal 
T}\left ( S\right )$. Again, as $\lambda _{n}$ 
converges to $\lambda $ for the Hausdorff 
topology, $\theta _{\lambda _{n}}$ and its 
tangent maps uniformly converge to $\theta 
_{\lambda }$ and its tangent maps as $n$ tends to 
$\infty $. 
\nobreak\par\nobreak 	We are now ready to prove 
the continuity property for $\Phi _{\rho _{0}}$. 
Let $\rho _{n}\in {\Cal R}\left ( S\right )$, 
$n\in {\Bbb N}$, be a sequence of representations 
converging to $\rho _{0}$. Let $\left ( 
m_{n},b_{n}\right )=\varphi ^{-1}\left ( \rho 
_{n}\right )\in {\Cal T}\left ( S\right )\times 
{\Cal M}{\Cal L}\left ( S\right )$, and let $\dot 
m_{n}\in T_{m_{n}}{\Cal T}\left ( S\right )$ be a 
sequence of tangent vectors converging to some 
$\dot m_{0}\in T_{m_{0}}{\Cal T}\left ( S\right 
)$. We want to show that $\Phi _{\rho _{n}}\left 
( \dot m_{n}\right )$ converges to $\Phi _{\rho 
_{0}}\left ( \dot m_{0}\right )$. 
\nobreak\par\nobreak 	For each $n$, let $\lambda 
_{n}$ be a maximal geodesic lamination containing 
the support of $b_{n}$. Extracting a subsequence 
if necessary, we can assume that $\lambda _{n}$ 
converges for the Hausdorff topology to some 
maximal geodesic lamination $\lambda _{0}$ 
containing the support of $b_{0}$. Let $U$ be a 
train track carrying $\lambda _{0}$. Then, by 
definition of all the maps involved,
$$\varphi _{\lambda _{n}}\left ( 
m_{n},b_{n}\right ) = \psi _{\lambda _{n}}\left ( 
\theta _{\lambda _{n}}^{-1}\left ( m_{n}\right 
)+iB_{n}\right )$$
for $n$ sufficiently large, where $B_{n}\in {\Cal 
W}\left ( U;{\Bbb R}/2\pi {\Bbb Z}\right )$ is 
the edge weight system corresponding to $b_{n}\in 
{\Cal H}\left ( \lambda _{n};{\Bbb R}/2\pi {\Bbb 
Z}\right )$. It follows that 
$$\Phi _{\rho _{n}}\left ( \dot m_{n}\right 
)=T_{\left ( m_{n},b_{n}\right )}\varphi \left ( 
\dot m_{n},0\right )=T_{\left ( m_{n},b_{n}\right 
)}\varphi _{\lambda _{n}}\left ( \dot 
m_{n},0\right )
	=T_{\left ( \theta _{\lambda _{n}}^{-1}\left ( 
m_{n}\right )+iB_{n}\right )} \psi _{\lambda 
_{n}}\left ( T_{m_{n}}\theta _{\lambda 
_{n}}^{-1}\left ( \dot m_{n}\right )\right ) .$$
By uniform convergence of the tangent maps, we 
conclude that $\Phi _{\rho _{n}}\left ( \dot 
m_{n}\right )$ converges to $\Phi _{\rho 
_{0}}\left ( \dot m_{0}\right )$ as $n$ tends to 
$\infty $. 
\nobreak\par\nobreak 	This concludes the proof of 
Lemma~13.  \qed\enddemo 
 
\nobreak\par\nobreak 	By Lemma~13, $\Phi _{\rho 
_{0}}$ depends continuously on $\rho _{0}$. In 
particular, its image $P_{\rho _{0}}$ depends 
continuously on $\rho _{0}$. Therefore the 
projection $\Pi _{\rho _{0}}:T_{\rho _{0}}{\Cal 
R}\left ( S\right )\rightarrow P_{\rho _{0}}$ 
parallel to $iP_{\rho _{0}}$ also depends 
continuously on $\rho _{0}$. This proves that the 
tangent map $T_{\rho _{0}}\left ( Q\circ \varphi 
^{-1}\right )=\Phi _{\rho _{0}}^{-1}\circ \Pi 
_{\rho _{0}}$ depends continuously on $\rho 
_{0}$, and concludes the proof of Proposition~12 
and Theorem~1.    \qed\enddemo

\head {\S }6. The map $\mu $ is not necessarily 
twice differentiable\endhead    
\nobreak\par\nobreak 	It is not difficult to show 
by explicit computations that the map $\mu $ is 
not necessarily twice differentiable. For 
instance, we can borrow such computations from 
\cite{PaS}. Let $S$ be a once punctured torus. On 
$S$, choose a hyperbolic metric $m_{0}\in {\Cal 
T}\left ( S\right )$ and a pair of simple closed 
$m_{0}$--geodesics $\gamma $, $\delta $ on $S$ 
meeting transversely in one point. If $\rho \in 
{\Cal R}\left ( S\right )$ is geometrically 
finite and if $M$ is the corresponding hyperbolic 
3--manifold, the boundary $\partial C_{M}$ is the 
union of two copies $\partial ^{+}C_{M}$ and 
$\partial ^{-}C_{M}$ of $S$, where the 
identification of $S$ with $\partial ^{+}C_{M}$ 
(resp. $\partial ^{-}C_{M}$) respects (resp. 
reverses) the orientation.  Let $\gamma ^{\pm }$ 
and $\delta ^{\pm }$ denote the closed geodesics 
of $\partial ^{\pm }C_{M}$ homotopic to $\gamma $ 
and $\delta $, respectively.
\nobreak\par\nobreak 	For $t\in {\Bbb R}$, let 
$\gamma _{t}\in {\Cal H}\left ( \lambda ;{\Bbb 
R}/2\pi {\Bbb Z}\right )$ be the Dirac transverse 
measure for $\gamma $ with mass the $\roman 
{mod}\thinspace 2\pi $ reduction of $t$, and let 
$\rho _{t}=\varphi _{\gamma }\left ( m_{0}, 
\gamma _{t}\right )$. The representation $\rho 
_{0}$ is Fuchsian, and defines a hyperbolic 
3--manifold $M_{0}$. For $t$ close to 0, we can 
then consider the hyperbolic metric $M_{t}\in 
{\Cal Q}{\Cal D}\left ( M_{0}\right )$ 
corresponding to $\rho _{t}$. 
\nobreak\par\nobreak 	First consider the case 
where $t$ is non-negative, and close to 0. Then, 
$\partial ^{+}C_{M_{t}}$ has induced metric 
$m_{0}$ and bending measured geodesic lamination 
$\gamma _{t}$. If we make the additional 
assumption that $\gamma $ and $\delta $ meet 
orthogonally for the metric $m_{0}$, it is shown 
in \cite{PaS} that $\partial ^{-}C_{M}$ is bent 
along $\delta ^{-}$; this can also be seen from 
symmetry considerations.
\nobreak\par\nobreak 	For $t\leqslant 0$ close to 
0,  it is now $\partial ^{-}C_{M_{t}}$ which has 
induced metric $m_{0}$ and bending measured 
lamination $\gamma _{-t}$, and $\partial 
^{+}C_{M}$ is bent along $\delta ^{+}$. In 
addition, the central equality of \cite{PaS} 
shows that the lengths of $\gamma ^{-}$ and 
$\delta ^{+}$ are related to $t$ by the formula
$$\cos^{2}\left ( {t/2}\right ) = 
\cosh^{2}l_{t}\left ( \gamma ^{-}\right ) 
\tanh^{2}l_{t}\left ( \delta ^{+}\right ). $$
Noting that $l_{t}\left ( \gamma ^{-}\right 
)=l_{0}\left ( \gamma \right )$, we conclude that 
$\tanh^{2}l_{t}\left ( \delta ^{+}\right 
)=\cos^{2}\left ( {t/2}\right 
)/\cosh^{2}l_{0}\left ( \gamma \right )$.
\nobreak\par\nobreak 	Therefore, for $t$ small, 
the function $\tanh^{2}l_{t}\left ( \delta 
^{+}\right )$ is equal to $\tanh^{2}l_{0}\left ( 
\delta \right )=1/\cosh^{2}l_{0}\left ( \gamma 
\right )$ if $t\geqslant 0$ and to $\cos^{2}\left 
( {t/2}\right )/\cosh^{2}l_{0}\left ( \gamma 
\right )$ if $t\leqslant 0$. This function of $t$ 
is not twice differentiable at 0. On the other 
hand, the curve $t\mapsto M_{t}$ is real analytic 
in ${\Cal Q}{\Cal D}\left ( M_{0}\right )$. It 
follows that $\mu $ is not twice differentiable 
at $M_{0}$.

\enddocument